\documentclass[12pt]{amsart}%
\usepackage[pagewise]{lineno}
\usepackage{graphicx,subcaption}
\usepackage{tikz}
\usepackage{tikz-cd}
\usetikzlibrary{matrix,arrows,decorations.pathmorphing}
\usepackage{amscd}
\usepackage{geometry}
\usepackage{amsmath}
\usepackage{amsfonts}
\usepackage{amssymb}%
\usepackage{enumerate}
\usepackage{yfonts}
\usepackage[shortlabels]{enumitem}
\usepackage{extpfeil}
\usepackage{cleveref}% http://ctan.org/pkg/cleveref
\providecommand{\U}[1]{\protect\rule{.1in}{.1in}}
\newtheorem{theorem}{Theorem}[section]
\theoremstyle{plain}

\newtheorem{corollary}[theorem]{Corollary}

\newtheorem{lemma}[theorem]{Lemma}
\newtheorem{question}{Question}
\newtheorem{proposition}[theorem]{Proposition}
\theoremstyle{definition}
\newtheorem{definition}{Definition}[section]

\newtheorem{remark}{Remark}

\numberwithin{equation}{section}

\newcommand{\step}[1]{\medskip\noindent\textit{#1.}}
\usepackage [autostyle, english = american]{csquotes}
\MakeOuterQuote{"}

\makeatletter
\def\oversortoftilde#1{\mathop{\vbox{\m@th\ialign{##\crcr\noalign{\kern3\p@}%
  \sortoftildefill\crcr\noalign{\kern3\p@\nointerlineskip}%
  $\hfil\displaystyle{#1}\hfil$\crcr}}}\limits}

\def\sortoftildefill{$\m@th \setbox\z@\hbox{$\braceld$}%
  \braceld\leaders\vrule \@height\ht\z@ \@depth\z@\hfill\braceru$}

\makeatother

\usepackage{scalerel,stackengine}
\stackMath
\newcommand\reallywidehat[1]{%
\savestack{\tmpbox}{\stretchto{%
  \scaleto{%
  \scalerel*[\widthof{\ensuremath{#1}}]{\kern-.6pt\bigwedge\kern-.6pt}%
  {\rule[-\textheight/2]{1ex}{\textheight}}%WIDTH-LIMITED BIG WEDGE
  }{\textheight}%
}{0.5ex}}%
\stackon[1pt]{#1}{\tmpbox}%
}
\begin{document}
\title[Nonembeddable contractible open manifolds]{Nonembeddable contractible open manifolds arising from Whitehead doubling}
\author{Shijie Gu}
\address{Department of Mathematics,
Northeastern University, Shenyang, Liaoning, China, 110004}
\email{shijiegutop@gmail.com}

\author{Jian Wang}
\address{State Key Laboratory of Mathematical Sciences,
Academy of Mathematics and Systems Science, Chinese Academy of Sciences,
Beijing 100190, China}
\email{jian.wang.4@amss.ac.cn}

\author{Yanqing Zou}
\address{School of Mathematical Sciences, Key Laboratory of MEA(Ministry of Education) \& Shanghai Key Laboratory of PMMP,
East China Normal University, Shanghai, China, 200241}
\email{yqzou@math.ecnu.edu.cn}
\thanks{}
\date{\today}

\begin{abstract}
We study a family of genus-one contractible open manifolds constructed by iterated Whitehead doubling from a nontrivial knot $K$ and an even half-twist $m$.  For each such pair, we prove that the resulting contractible open $3$-manifold $W(K,m)$ does not embed as an open subset of any compact, locally connected and locally $1$-connected metric $3$-space.  We also classify this family: $W(K,m)$ is homeomorphic to $W(K',m')$ if and only if $m=m'$ and $K$ is isotopic to $K'$.  Thus the knot type and the half-twist parameter form a complete invariant for these nonembeddable contractible open manifolds.  The proof combines the topology of ends with JSJ decompositions, hyperbolic pieces, and a rank estimate for iterated Whitehead doubled knot groups.  The same method gives infinitely many pairwise non-homeomorphic higher-dimensional examples which embed in no compact, locally connected and locally $1$-connected metric space of the same dimension.
\end{abstract}

\maketitle

\section{Introduction}
The study of contractible open $3$-manifolds begins with Whitehead's construction of a contractible open $3$-manifold not homeomorphic to $\mathbb R^3$ \cite{Whi35}.  Whitehead's example showed that contractibility alone does not control the topology at infinity.  Subsequent work of McMillan produced large families of such manifolds, in particular genus-one examples obtained as increasing unions of nested solid tori \cite{McM62}.  Since then, contractible open $3$-manifolds have served as test objects for several basic questions in low-dimensional topology, including embeddability in compact spaces, covering-space behavior, and the extent to which defining sequences determine the resulting open manifold; see, for example, \cite{Mye88, Mye99, Wri92}.

This paper is concerned with two of these questions. The first is an embeddability problem: when can a contractible open $3$-manifold occur as an open subset of a compact space?  A classical, more restrictive version is the missing-boundary problem, in which a
manifold is obtained from a compact manifold by deleting a closed subset of its boundary.  Tucker characterized $P^2$-irreducible
missing-boundary $3$-manifolds in terms of finite generation of the
fundamental groups of complements of compact polyhedra
\cite{Tuc74}.  Simon conjectured that every covering of a compact
irreducible $3$-manifold with finitely generated fundamental group
is missing-boundary and proved this in several cases \cite{Sim76}.
The conjecture is now a consequence of the Tameness and
Geometrization Theorems.  The missing-boundary property is stronger than arbitrary open embeddability in a compact $3$-manifold.  Kister and McMillan exhibited locally Euclidean factors of $\mathbb E^4$ that cannot be embedded in $\mathbb S^3$, and Haken's finiteness theorem obstructs embedding many such examples in compact $3$-manifolds \cite{KM62, Hak68}.  For open $3$-manifolds exhausted by compact submanifolds with torus boundary, Messer and Wright gave necessary and sufficient geometric conditions for such embeddings \cite{MW79}.  Sternfeld constructed examples that do not embed in the broader category of compact, locally connected and locally $1$-connected metric spaces \cite{Ste77}; more recently, Gu proved the same for Bing's contractible open manifold and related examples \cite{Gu21}.

% The first is an embeddability problem: when can a contractible open $3$-manifold occur as an open subset of a compact space?  Kister and McMillan exhibited locally Euclidean factors of $\mathbb E^4$ that cannot be embedded in $\mathbb S^3$; one such example is Bing's manifold (see Figure~\ref{Bing's manifold}) \cite{KM62}. Haken's finiteness theorem provides a classical obstruction to embedding many such examples in compact $3$-manifolds \cite{Hak68}. Messer and Wright studied the general problem of embedding open $3$-manifolds in compact $3$-manifolds \cite{MW79}.  A natural question asks whether such manifolds embed in compact, locally connected and locally $1$-connected metric spaces.  Sternfeld constructed examples (see Figure \ref{Sternfeld's manifold}) that do not embed in this broader category \cite{Ste77}; more recently, Gu used Sternfeld's methods to prove nonembeddability for Bing's contractible open manifold and related examples \cite{Gu21}. 
%In a different compact three-dimensional setting, Belletti--Detcherry recently
% used quantum invariants to obstruct embeddings between compact oriented
% $3$-manifolds \cite{BD26}.  For the manifolds constructed below, however, every compact $3$-submanifold lies in a solid-torus stage and hence embeds in $S^3$. Their nonembeddability arises only from the infinite nesting of these stages and is therefore encoded in the topology at infinity.

There is also a parallel smooth four-dimensional line of work in which Floer- and gauge-theoretic invariants obstruct embeddings of closed $3$-manifolds.  For example, Taniguchi used instanton Floer theory to obstruct certain embeddings of homology $3$-spheres into homology $S^3\times S^1$ \cite{Tan18}, while Daemi--Lidman--Eismeier constructed infinitely many closed $3$-manifolds which do not embed in any closed symplectic $4$-manifold, using Heegaard Floer correction terms and instanton moduli spaces \cite{DLE24}; see also \cite{ET23, AMP24}.  In a different same-dimensional compact direction, Belletti--Detcherry recently used quantum invariants to obstruct embeddings between compact oriented $3$-manifolds \cite{BD26}.  These results concern smooth $4$-manifolds or compact $3$-manifolds, whereas the present paper treats contractible open $3$-manifolds embedded in compact $3$-dimensional metric spaces and uses rank growth and JSJ structure at infinity.

We first describe the family.  Let $K$ be a nontrivial knot and let $m$ be an even integer.  Throughout the construction, $K$ is represented in the local cube by a fixed diagram $D_K$, and the normal framing of this local model is its blackboard framing.  The integer $m$ denotes the number of half-twists inserted in the box of Figure~\ref{3_1knot}, measured relative to this blackboard framing.  Thus the ordinary satellite twisting number of the corresponding Whitehead double is
\[
   \tau=\frac m2+\operatorname{wr}(D_K),
\]
where $\operatorname{wr}(D_K)$ is the writhe of the fixed local diagram.  The notation $W(K,m)$ always refers to the contractible open $3$-manifold constructed from this blackboard-framed local model by the direct-limit construction of Section~\ref{Section: nonembeddability}.  At each stage a solid torus is inserted by the same Whitehead-pattern layer, with the knotted part modeled on $K$.

The first main result is a nonembeddability theorem for the whole family.

\begin{theorem}\label{Thm: nonembeddability}
Let $W^3=W(K,m)$ be a contractible open manifold constructed as in Section~\ref{Section: nonembeddability}, where $K$ is a nontrivial knot and $m$ is an even integer. Then $W^3$ does not embed as an open subset of any compact, locally connected, locally $1$-connected metric $3$-space. In particular, $W^3$ embeds in no compact $3$-manifold.
\end{theorem}

This theorem produces a large family of contractible open 3-manifolds that embed in no compact, locally connected, locally 1-connected
metric 3-space, in analogy with Haken’s finiteness theorem.  Kister and McMillan asked whether certain exotic contractible open manifolds embed in compact $3$-manifolds.  Haken's finiteness theorem \cite[P. 65--69]{Hak68} gives a broad nonembeddability result for examples of the type shown in Figure~\ref{3_1knot}, resolving their conjecture \cite{KM62}: a compact $3$-manifold contains only finitely many nonparallel incompressible surfaces.  It is more challenging to prove nonembeddability in general compact metric spaces which are locally connected and locally $1$-connected, since Haken's finiteness theorem is no longer available in that setting.  Instead, the proof of Theorem~\ref{Thm: nonembeddability} uses rank growth in knot groups associated to the end of $W(K,m)$.  In particular, replacing the trefoil-knotted hole in Figure~\ref{3_1knot} with a $K$-knotted hole yields nonembeddable examples, thereby resolving \cite[Question~2]{Gu21}. Although the rank-growth argument also shows that $W(K,m)$ fails Tucker's missing-boundary criterion, this alone does not rule out an arbitrary open embedding in a compact $3$-manifold.  The
stronger conclusion uses the ambient compact space: such an embedding
would imply that the complement of a fixed initial solid-torus stage
has finitely generated fundamental group and surjects onto a sequence
of knot groups whose ranks tend to infinity. 

The second question is a classification problem.  Many constructions of genus-one contractible open $3$-manifolds are organized by defining sequences of nested solid tori.  Geometric index is a basic invariant in this setting, going back to Schubert's work on knots in solid tori and used in McMillan-type constructions.  It is also central in later work on Bing--Whitehead Cantor sets and generalized Whitehead links \cite{GRWZ11, GRW18}.  In this paper, however, all examples have the same geometric-index sequence: each consecutive inclusion of defining solid tori has geometric index \(2\).  Thus the classification
cannot rely on geometric index alone.  Instead, it recovers finer data at infinity: the knot type
used to knot the hole and the twisting number in the Whitehead-pattern layer.

\begin{theorem}[Classification]\label{Thm: homeom of manifolds intro}
Let $K$ and $K'$ be two nontrivial knots and let $m$ and $m'$ be even integers, with $W(K,m)$ and $W(K',m')$ constructed as in Section~\ref{Section: nonembeddability}. Then $W(K, m)$ is homeomorphic to $W(K',m')$ if and only if
\begin{equation*}
  m=m' \quad \text{ and }\quad  K \text{ is isotopic to } K'.
\end{equation*}
\end{theorem}

Theorem~\ref{Thm: homeom of manifolds intro} is another main contribution of the paper.  It says that the end of $W(K,m)$ remembers the complete pair $(K,m)$ (not merely coarse data such as geometric index).  The proof uses JSJ decompositions of finite layers at infinity.  Each layer contains a distinguished, or marked, hyperbolic component.  The marked components are intrinsic to the end; their complementary JSJ pieces recover the knot complement associated to $K$, and hence recover $K$ by the knot complement theorem.  The parameter $m$ is recovered from boundary-slope data on the marked Whitehead-pattern component.  More concretely, the preferred longitude on the outer cusp and the Seifert-fiber slope on the adjacent connected-sum component determine a Dehn filling which is the exterior of the $m/2$-twist knot.  Thus the pair $(K,m)$ is a complete invariant for the homeomorphism type within this family.

This classification should be compared with the double $3$-space problem for contractible open $3$-manifolds.  Gabai proved that the Whitehead manifold is the union of two copies of $\mathbb R^3$ whose intersection is again homeomorphic to $\mathbb R^3$ \cite{Gab11}.  Garity--Repov\v{s}--Wright then used generalized Whitehead links, geometric index, and interlacing theory to construct uncountably many contractible $3$-manifolds with the double $3$-space property and uncountably many without it \cite{GRW18}.  Their work shows that geometric index and interlacing provide powerful tools for distinguishing large classes of genus-one contractible open $3$-manifolds.  Our problem and methods are different: rather than deciding the double $3$-space property, we classify a family with fixed geometric index by recovering knot-theoretic and hyperbolic-geometric data from the end.

We now recall the algebraic ingredient used in the proof of Theorem~\ref{Thm: nonembeddability}.  The \emph{rank} of a group $G$, denoted $r(G)$, is the smallest cardinality of a generating set of $G$.  For a knot $K\subset S^3$, let $G(K)$ be its knot group and denote the rank of $G(K)$ by $r(K)$.  In \cite{Gu21}, the first author proved that Bing's manifold (see \ref{Bing's manifold}) and Sternfeld's manifold (see \ref{Sternfeld's manifold}) do not embed in any compact, locally connected and locally $1$-connected metric space.  Those examples are built from the trefoil-knotted hole of Figure~\ref{3_1knot}, and the proof depends on rank growth for iterated Whitehead doubles of the trefoil.  That growth was obtained by covering-space methods and computer-assisted calculations.

Extending this approach to an arbitrary $K$-knotted hole leads naturally to difficult quotient questions for knot groups, for instance the following.
\begin{question}
Does there exist a nonabelian knot group $G$ such that $\mathbb{A}_n$ is not a quotient of $G$ for every $n>3$?
\end{question}
A full answer appears beyond reach; see \cite{Gu17} for discussion.  See also \cite{BBK21, BKM24} for results using symmetric and Coxeter group quotients to bound the meridional rank of knot groups.  Hence it is unclear how to extend the methods of \cite{Gu21} to genus-one contractible open manifolds constructed from an arbitrary $K$-knotted hole.

The rank estimate below provides a different route.  It replaces the trefoil-specific computation in \cite{Gu21} with a uniform JSJ-theoretic argument.

\begin{proposition}\label{Thm: Lower bound for WD knot group}
  Let $K$ be a nontrivial knot and $n\in \mathbb{N}$. Then $r\left(\operatorname{WD}^n(K)\right)\geq n+1$.
\end{proposition}

The complement of a tubular neighborhood of $\operatorname{WD}^n(K)$ admits a JSJ decomposition with at least $n$ hyperbolic pieces coming from the Whitehead-doubling layers.  Weidmann's theorem \cite[Thm. 5]{Wei02}, in the form stated in Lemma~\ref{Lemma: lower bound for the rank}, turns these hyperbolic pieces into the lower bound in Proposition~\ref{Thm: Lower bound for WD knot group}.  This proposition supplies the missing piece of the puzzle needed to adapt the contradiction argument of \cite{Gu21} to the present setting and hence to prove Theorem~\ref{Thm: nonembeddability}.

The same method also gives higher-dimensional examples.

\begin{theorem}\label{Thm: high dimensional nonembeddability}
For every $n\geq 3$ there exist infinitely many non-homeomorphic contractible open $n$-manifolds $W^n$ which embed in no compact, locally connected, locally $1$-connected metric $n$-space. In particular, $W^n$ embeds in no compact $n$-manifold.
\end{theorem}

The paper is organized as follows.  Section~\ref{Section: Proof of the main theorem} proves Proposition~\ref{Thm: Lower bound for WD knot group} and records consequences for tunnel number, asymptotic rank, and the rank-versus-genus problem for a large subfamily.  Section~\ref{Section: nonembeddability} constructs the manifolds $W(K,m)$ and proves Theorem~\ref{Thm: nonembeddability}.  The same section proves the classification theorem by recovering the knot type from marked hyperbolic components and the half-twist parameter from boundary slopes.  It concludes with the higher-dimensional extension and further examples obtained by inserting Whitehead-manifold layers.

\step{Acknowledgments} The authors would like to thank Matthew Hedden, Sebastian Baader, Hongbin Sun, Jiming Ma and Louis Funar for helpful conversations and encouragement. The first author was supported in part by NSFC grant 12201102. The second author thanks the Simons Center for Geometry and Physics for its hospitality, where a part of this work was completed during his visit.  The third author was supported by NSFC 12326601, 12471065 and partially by Science and Technology Commission of Shanghai Municipality (STCSM), grant No. 22DZ2229014.
\section{Rank estimates for iterated Whitehead doubles}
\label{Section: Proof of the main theorem}

 \begin{definition}\label{Def: Whitehead doubling}
 Let $K_P\subset V_P\subset S^3$ be a knot contained in an unknotted solid torus $V_P\subset S^3$, such that $K_P$ is not contained in any $3$-ball in $V_P$. The pair $(V_P, K_P)$ is called a \emph{pattern} and $K_P$ is referred to as the \emph{pattern knot}.

Let $K_C \subset S^3$ be a knot, and $V_C$ be a tubular neighborhood of $K_C$ in $S^3$. Let
$h: V_P \to V_C$
be a homeomorphism onto $V_C$. The image
$K_W := h(K_P) \subset V_C \subset S^3$
is called a \emph{satellite knot} with \emph{companion knot} $K_C$ and pattern $(V_P, K_P)$.

\end{definition}

\begin{figure}[h!]
  \centering
\begin{subfigure}{0.4\textwidth}
  \centering
  \includegraphics[width=0.9\textwidth]{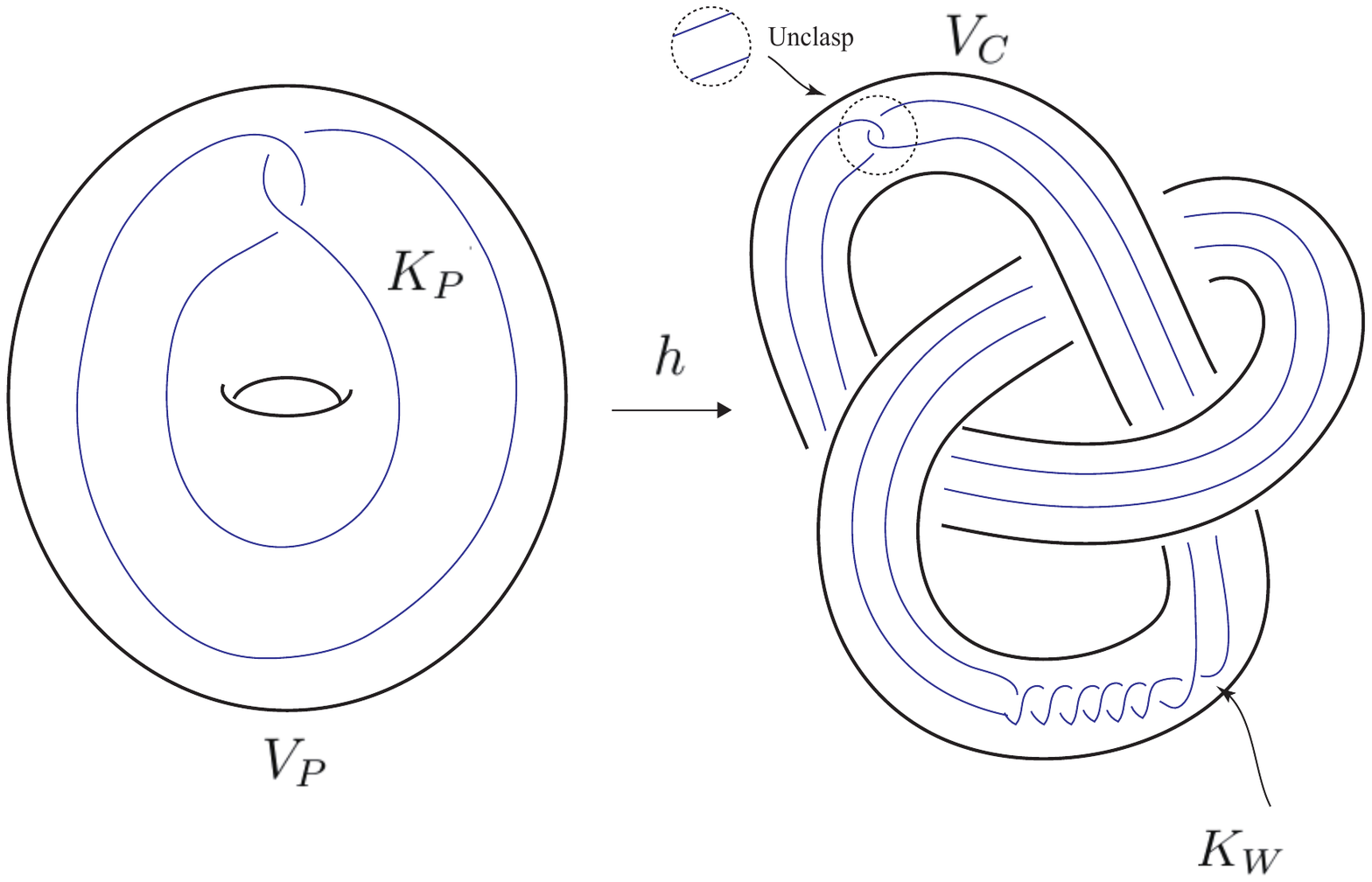}
  \caption{An untwisted Whitehead double of a trefoil knot.}
  \label{whiteheaddouble_3_1_untwisted}
\end{subfigure}\hspace{2em}
\begin{subfigure}{0.4\textwidth}
  \centering
  \includegraphics[width=0.9\textwidth]{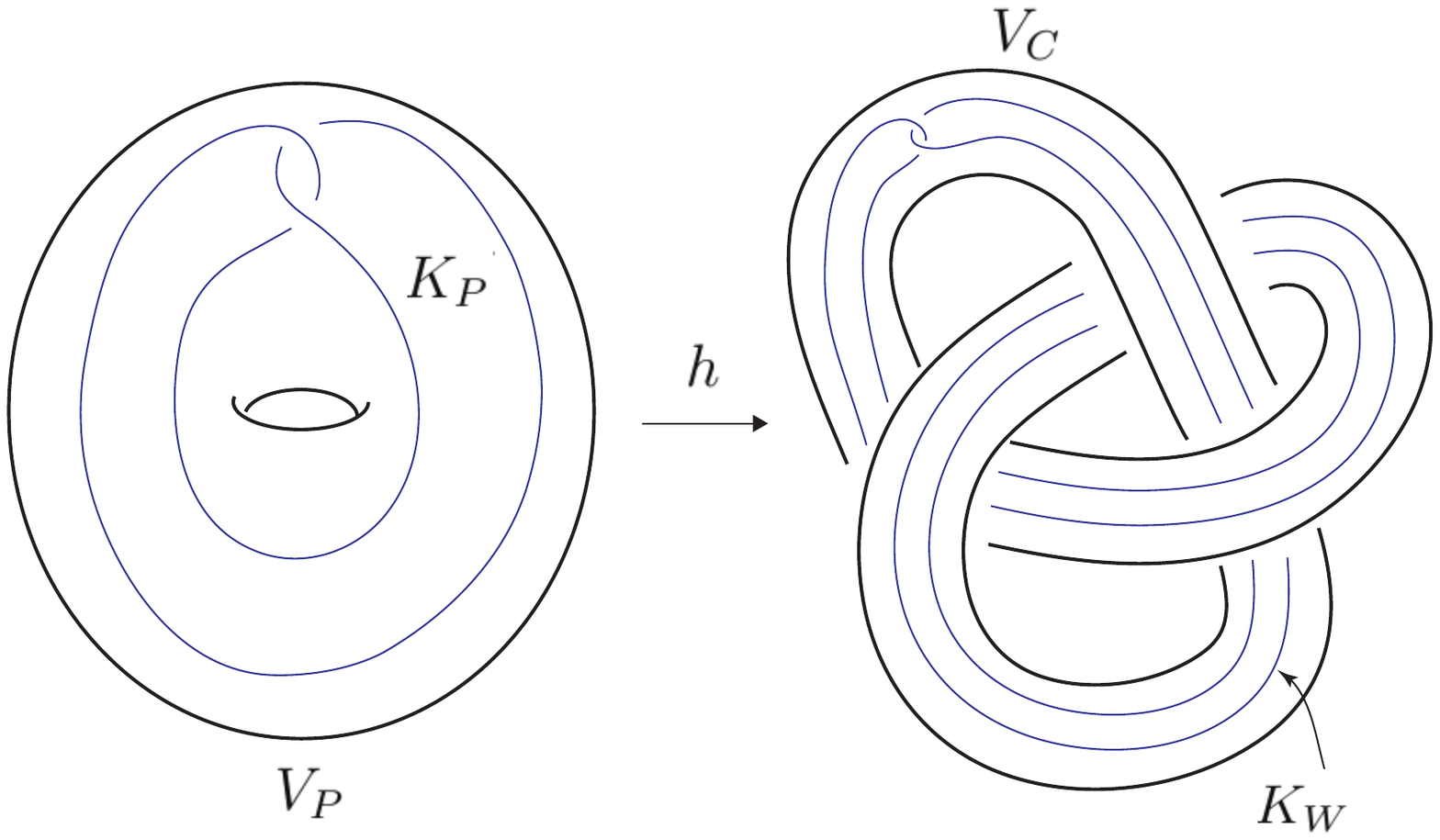}
  \caption{A 3-twisted Whitehead double of a trefoil knot.}
  \label{whiteheaddouble_3_1}
\end{subfigure}
\caption{A knot $K_W$ with a trefoil knot as companion.}
\label{A knot K_W with a trefoil knot as companion}
\end{figure}

Different choices of the homeomorphism $h$ may result in different knots $K_W$. The knot type of $K_W$ is determined by the twisting number, which encodes how $h$ maps the longitude and meridian of $V_P$ into those of $V_C$. Unclasping and reconnection for $K_W$ shown in Figure \ref{whiteheaddouble_3_1_untwisted} produces  a two-component link that forms the boundary of a closed ribbon. The linking number is referred to as the \emph{twisting number} \cite[P. 166]{Rol76}.

\begin{definition}\label{whitehead-double-def} Let  $(K_P, V_P)$ be the pattern as illustrated  in Figure \ref{whiteheaddouble_3_1_untwisted} and $K_C\subset S^3$ be any knot. We say that $K_{W}$ is an \emph{untwisted Whitehead double} of $K_{C}$ if the twisting number is zero.\footnote{In this case, the homeomorphism $h: V_P \to V_C$ is faithful, meaning that $h$ takes the preferred longitude and meridian of $V_P$ respectively to the preferred longitude and meridian of $V_C$.}  Otherwise, we call $K_W$ a \emph{twisted Whitehead double} of $K_C$.
\end{definition}
For instance, the satellite knot in Figure \ref{whiteheaddouble_3_1} is a 3-twisted Whitehead double of a trefoil knot, where the twisting arises from the writhe of a trefoil knot.

\begin{remark}In what follows, we do not specify whether a (possibly iterated) Whitehead doubling is twisted, since Proposition~\ref{Thm: Lower bound for WD knot group} applies in both the twisted and untwisted cases. Indeed, the JSJ pieces of the resulting doubled knots remain hyperbolic regardless of the twisting. See also the proof of Lemma~\ref{Lemma: toral differences are hyperbolic}.

\end{remark}

For a nontrivial knot $K$ in $S^3$, we consider its $n$th iterated Whitehead doubling, denoted by $\operatorname{WD}^n(K)$. We use this construction to identify the JSJ decomposition of the complement and, in particular, its hyperbolic pieces.  For this purpose, we introduce a nested sequence of solid tori $T_{1} \supset T_{2} \supset \cdots \supset T_{n}$ whose boundary tori provide the required incompressible surfaces.

Recall that $(K_P, V_{P})$ denotes the pattern pair illustrated in Figure \ref{whiteheaddouble_3_1_untwisted} and there exists a homeomorphism $h_1$:
\begin{equation*}
h_1: (K_P, V_P)\rightarrow (\operatorname{WD}(K), T_1),
\end{equation*} where $T_1$ ($= V_C$) is the tubular neighborhood of $K$ in $S^3$. Choosing a tubular neighborhood $T_2$ of $\operatorname{WD}(K)$ ($=K_W$) in $\operatorname{Int}T_1$, we obtain another homeomorphism $h_2:  (K_P, V_P)\rightarrow (\operatorname{WD}^2(K), T_2) $. Iterating this process, we construct a nested sequence $\{T_{i}\}_{i=1}^n$ such that $T_i$ is a tubular neighborhood of $\operatorname{WD}^{i-1}(K)$.

Since the original knot $K$ is nontrivial, each $\operatorname{WD}^{i}(K)$ is nontrivial in $S^3$. Combining this observation with the Seifert--van Kampen theorem, we will show that $\partial T_i$ is incompressible in $S^3\setminus\operatorname{WD}^n(K)$ for all $i\leq n$.

\begin{lemma}\label{Lemma: incompressible torus boundary}
  Let $K$ be a nontrivial knot and  $\{T_{i}\}_{i=1}^{n}$ be constructed as above. Then, for any $i\leq n$, the torus $\partial T_i$  is incompressible in $S^{3} \setminus  \operatorname{WD}^{n}(K)$.
\end{lemma}
\begin{proof} Recall that $h_i: (K_P, V_P)\rightarrow (\operatorname{WD}^i(K), T_i)$ is a homeomorphism for each $i$. Consequently, for $i\leq n-1$, the region $T_i\setminus \operatorname{Int}T_{i+1}$ is homeomorphic to the exterior of the Whitehead pattern in $V_P$. Since the interior of this pattern exterior is the complement of a Whitehead link with twisting in $S^3$, it follows that both $\partial T_i$ and $\partial T_{i+1}$ are incompressible in $T_i\setminus \operatorname{Int}T_{i+1}$. Moreover, we observe that the inclusion \begin{equation*}
  \pi_1(\partial T_1)\rightarrow  \pi_1(S^3\setminus \operatorname{Int}T_1)
\end{equation*} is injective, since $K$ is knotted in $S^3$.

By Theorem 11.60 of \cite{Rot} (a version of the Seifert--van Kampen theorem), we obtain the following amalgamated  free product decomposition:
\begin{equation*}
  \pi_1(S^3\setminus \operatorname{Int}T_2)\cong\pi_1(S^3\setminus \operatorname{Int}T_1)\ast_{\pi_1(\partial T_1)}\pi_1(T_1\setminus \operatorname{Int}T_2).
\end{equation*} In particular, the induced homomorphisms
\begin{equation*}
  \pi_1(T_1\setminus \operatorname{Int}T_2)\rightarrow \pi_1(S^3\setminus \operatorname{Int}T_2) \quad \text{and}\quad \pi_1(S^3\setminus \operatorname{Int}T_1)\rightarrow \pi_1(S^3\setminus \operatorname{Int}T_2)
\end{equation*} are injective, which implies that $\partial T_1$ and $\partial T_2$ are incompressible in $S^3\setminus \operatorname{Int}T_2$. Repeating this argument, we can conclude that $\partial T_i$ is incompressible in $S^3\setminus \operatorname{Int}T_n$ for all $i\leq n$.

Moreover, the homeomorphism $h_n: (K_P, V_P)\rightarrow (\operatorname{WD}^n(K), T_n)$ ensures that $\partial T_n$ is incompressible in $T_n\setminus \operatorname{WD}^n(K)$. The analogous amalgamated-product argument shows that the inclusion $\pi_1(S^3\setminus \operatorname{Int}T_n)\rightarrow \pi_1(S^3\setminus \operatorname{WD}^n(K))$ is injective, which yields the desired result.
\end{proof}

\begin{lemma}\label{Lemma: toral differences are hyperbolic}
  Let $K$ be a nontrivial knot and  $\{T_{i}\}_{i=1}^{n}$ be constructed as above. Then for any $i\leq n-1$, one has that $\operatorname{Int}T_{i} \setminus  {T}_{i+1}$ is hyperbolic in $S^{3} \setminus \operatorname{WD}^{n}(K)$.
\end{lemma}
\begin{proof} Recall that $h_i: (K_P, V_P)\rightarrow (\operatorname{WD}^i(K), T_i)$ is a homeomorphism.  Thus $\operatorname{Int}T_i\setminus T_{i+1}$ is homeomorphic to the interior of the Whitehead-pattern exterior in $V_P$.

By Thurston's hyperbolization theorem for the Whitehead-link exterior \cite{Thu}, $\operatorname{Int}T_i\setminus T_{i+1}$ admits a complete finite-volume hyperbolic structure.  By Mostow--Prasad rigidity \cite{Mos73, Pra73}, this determines a hyperbolic JSJ piece.
\end{proof}

\begin{figure}[h!]
  \centering
  \includegraphics[ width=12cm, height=7cm]{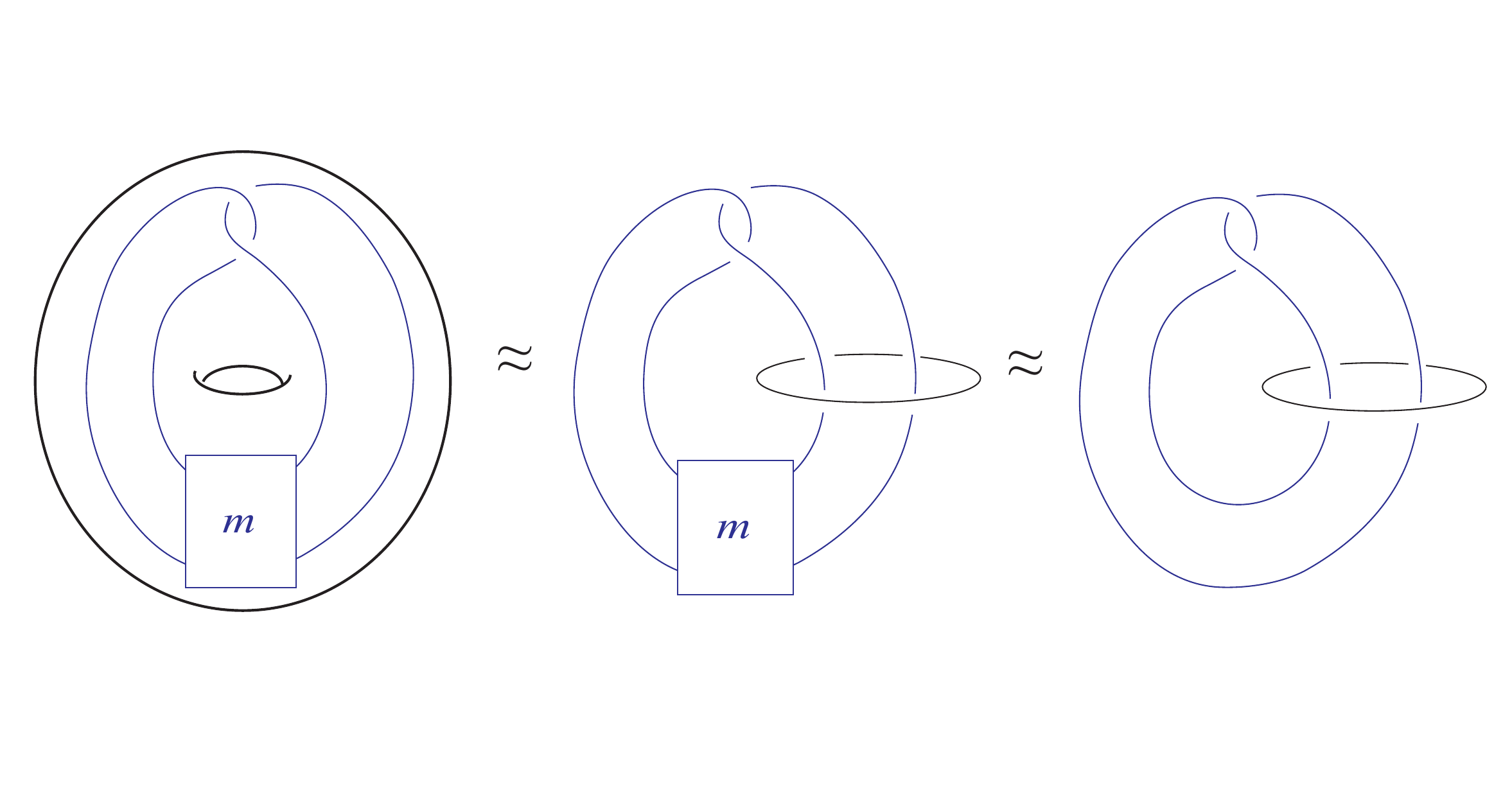}
  \caption{The Whitehead link complement in $S^3$ is homeomorphic to the complement of a twist knot with $m$ half-twists in a solid torus, where $m \in 2\mathbb{Z}$.}
  \label{twisted whitehead double}
\end{figure}

\begin{remark}
The left of Figure~\ref{twisted whitehead double} illustrates the embedding of $\operatorname{WD}^i(K)$ relative to $T_i$, where the blue curve represents the knot $\operatorname{WD}^{i}(K)$. This means that $\operatorname{Int}T_{i}\setminus T_{i+1}$ is the complement of the Whitehead link with half-twists. The ``twist number'' $m$ in Figure \ref{twisted whitehead double} is not the same as the ``twisting number'' of a twisted Whitehead doubling. The former counts the number of half-twists in a single component, whereas the latter is the linking number. However, by changing the number of full twists in the blue component of the link, one can adjust the twisting number of the corresponding Whitehead doubling. See Section \ref{Section: nonembeddability}.
\end{remark}

In light of Weidmann~\cite{Wei02}, we will exploit the hyperbolic pieces in $S^3 \setminus \mathrm{WD}^n(K)$ to obtain an estimate for the rank of the knot group of $\mathrm{WD}^n(K)$.

\begin{lemma}\label{Lemma: lower bound for the rank}
  Let $M^{3}$ be a compact orientable $3$-manifold with incompressible torus boundary, and let $n$ be the number of hyperbolic pieces in its JSJ decomposition.
  Then $r\left(\pi_1(M^3)\right) \geq n+1$.
\end{lemma}
\begin{proof}
  By performing suitable Dehn fillings on all torus boundary components, $M^{3}$ becomes a closed orientable $3$-manifold. The argument of \cite[Thm.5]{Wei02} then applies verbatim.
\end{proof}

We will first apply Lemma~\ref{Lemma: incompressible torus boundary} and Lemma~\ref{Lemma: toral differences are hyperbolic} to determine the JSJ decomposition of $S^3 \setminus \mathrm{WD}^n(K)$. We then invoke Lemma~\ref{Lemma: lower bound for the rank} to complete the proof of Proposition~\ref{Thm: Lower bound for WD knot group}.

\begin{proof}[Proof of Proposition \ref{Thm: Lower bound for WD knot group}] For a nontrivial knot $K$, consider a nested sequence of solid tori $\{T_i\}_{i=1}^n$, where $T_i$ is the tubular neighborhood of $\operatorname{WD}^{i-1}(K)$. To analyze  the JSJ decomposition of $S^3\setminus \operatorname{WD}^n(K)$, we  begin with  a JSJ decomposition of $S^3\setminus T_1$ and let $\Sigma_{k}$ denote the corresponding collection of incompressible surfaces.

 As shown in the proof of Lemma \ref{Lemma: incompressible torus boundary}, the inclusion  $S^3\setminus T_1\rightarrow S^3\setminus \operatorname{WD}^n(K)$ induces an injective homomorphism of fundamental groups:
$$\pi_1(S^3\setminus T_1)\rightarrow \pi_1(S^3\setminus \operatorname{WD}^n(K)).$$Hence,  the surfaces $\{\Sigma_k\}$ remain  incompressible in  $S^3\setminus \operatorname{WD}^n (K)$. Moreover, Lemma \ref{Lemma: incompressible torus boundary} implies  the boundary  $\partial T_i$ is incompressible in $S^3\setminus \operatorname{WD}^n (K)$ while  Lemma \ref{Lemma: toral differences are hyperbolic} shows that  the pieces $\operatorname{Int}T_i\setminus T_{i+1}$ are hyperbolic.

Combining these hyperbolic pieces  with the incompressible surface $\Sigma_k$ produces a JSJ decomposition of $S^3\setminus \operatorname{WD}^{n}(K)$ containing at least $n$ hyperbolic components. Applying  Lemma \ref{Lemma: lower bound for the rank} then establishes the required lower bound on the rank.
\end{proof}

Using the relation between the tunnel number and the rank, we obtain the following immediate consequence.

\begin{corollary}\label{Cor: tunnel number lower bound}
  Let $K$ be a nontrivial knot and $n\in \mathbb{N}$. Then the tunnel number $t\left(\operatorname{WD}^n(K)\right)\geq n$.
\end{corollary}
\begin{proof}
Let $E_n=S^3\setminus \operatorname{Int}N(\operatorname{WD}^n(K))$, and let $g(E_n)$ denote its Heegaard genus.  By definition,
\[
 t(\operatorname{WD}^n(K))=g(E_n)-1.
\]
A genus-$g(E_n)$ Heegaard splitting of $E_n$ gives a generating set of size $g(E_n)$ for $\pi_1(E_n)$.  Hence
\[
 r(\operatorname{WD}^n(K))\leq g(E_n)=t(\operatorname{WD}^n(K))+1.
\]
Proposition~\ref{Thm: Lower bound for WD knot group} gives $r(\operatorname{WD}^n(K))\geq n+1$, and therefore $t(\operatorname{WD}^n(K))\geq n$.
\end{proof}

We now examine the relationship between Whitehead doubling and the tunnel number, from which the asymptotic behavior of the rank follows.

\begin{corollary}\label{Cor: Linear growth of rank}
  Let $K$ be a nontrivial knot and $n\in \mathbb{N}$. Then one has that
  $$\lim_{n\rightarrow \infty} \frac{r\left(\operatorname{WD}^n(K)\right)}{n+1}=1.$$
\end{corollary}
\begin{proof}
For any knot $K'$, the tunnel number satisfies
\[
 t(\operatorname{WD}(K'))\leq t(K')+1.
\]
Iterating gives $t(\operatorname{WD}^n(K))\leq t(K)+n$.  Since the rank of a knot group is bounded above by the Heegaard genus of its exterior, we have
\[
 r(\operatorname{WD}^n(K))\leq t(\operatorname{WD}^n(K))+1\leq t(K)+n+1.
\]
Together with Proposition~\ref{Thm: Lower bound for WD knot group}, this gives
\[
 n+1\leq r(\operatorname{WD}^n(K))\leq n+1+t(K),
\]
and the stated limit follows.
\end{proof}

The following result provides numerous knots satisfying the famous rank versus genus problem.
\begin{corollary}\label{Cor: rank genus equality}
Let $K$ be a hyperbolic knot with tunnel number one. Then, for every $n\geq 1$,
\[
 r(\operatorname{WD}^n(K))=g(E(\operatorname{WD}^n(K)))=n+2,
 \qquad
 t(\operatorname{WD}^n(K))=n+1,
\]
where $E(\operatorname{WD}^n(K))$ denotes the exterior of $\operatorname{WD}^n(K)$.
\end{corollary}

\begin{proof}
The complement $S^3\setminus \operatorname{Int}T_1$ is hyperbolic, and the proof of Proposition~\ref{Thm: Lower bound for WD knot group} gives $n$ further hyperbolic pieces coming from the Whitehead-doubling layers.  Weidmann's theorem therefore gives
\[
 r(\operatorname{WD}^n(K))\geq n+2.
\]
Since $t(K)=1$, the tunnel-number estimate above gives
\[
 g(E(\operatorname{WD}^n(K)))=t(\operatorname{WD}^n(K))+1\leq n+2.
\]
The general inequality $r(\operatorname{WD}^n(K))\leq g(E(\operatorname{WD}^n(K)))$ forces equality throughout.
\end{proof}

\section{Nonembeddability and classification of the associated contractible open manifolds}\label{Section: nonembeddability}

In this section, we construct contractible $3$-manifolds via Whitehead doubling with respect to a knot $K$ and an even blackboard-framed half-twist parameter $m$. Using Proposition~\ref{Thm: Lower bound for WD knot group}, we show that such manifolds cannot be embedded into any compact, locally connected and locally $1$-connected metric $3$-space. Furthermore, varying either the knot type or this blackboard-framed half-twist parameter produces infinitely many distinct contractible $3$-manifolds, which may be distinguished by marked hyperbolic JSJ components and their boundary-slope data.

\subsection{General construction}
\begin{figure}[h!]
  \centering
  \includegraphics[ width=8cm, height=10cm]{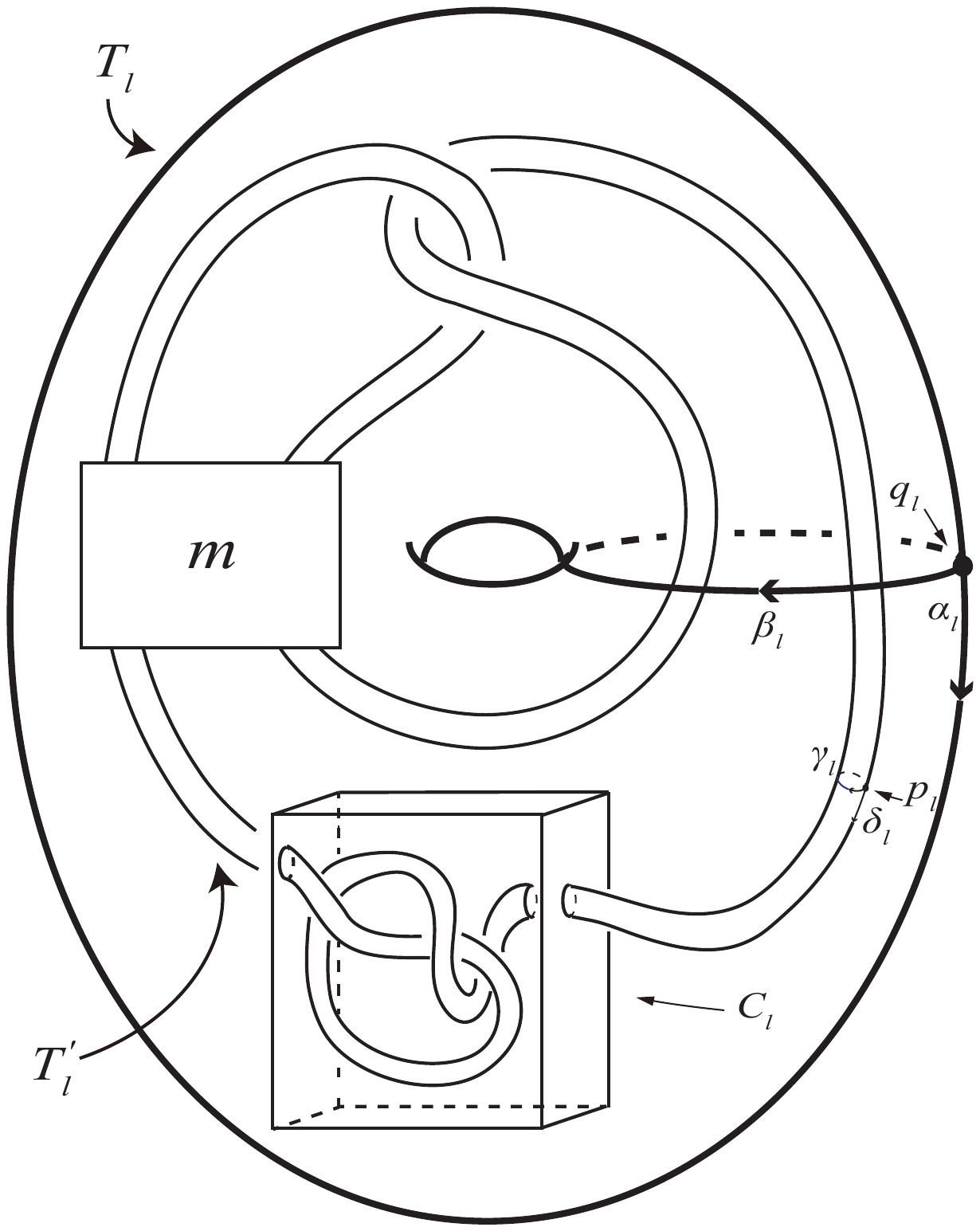}
  \caption{$L_l = T_l \setminus \operatorname{Int} T_l' = T_{l}^{\ast}\setminus \operatorname{Int} T_{l-1}^{\ast}$. The ``inner'' boundary component of $L_l$ is $\partial T_l'$. The ``outer'' boundary component
  of $L_l$ is $\partial T_l$. The box $m$ represents $m$ half-twists relative to the blackboard framing of the local $K$-diagram.}
  \label{3_1knot}
\end{figure}

\paragraph{Framing convention.}
We first fix the convention used for the integer $m$.  In each local cube $C_l$ appearing in Figure~\ref{3_1knot}, the knot $K$ is represented by a fixed diagram $D_K$, and the framing of the $K$-summand is the blackboard framing determined by this diagram.  The box labelled $m$ records $m$ half-twists relative to this blackboard framing; throughout the paper $m$ is required to be even.  We write $\operatorname{wr}(D_K)$ for the writhe of the fixed diagram.  With this convention, the actual twisting number of the Whitehead double obtained after the $K$-summand is inserted is
\[
  \tau= m/2+\operatorname{wr}(D_K).
\]
Thus $m$ is a diagrammatic, blackboard-framed layer parameter, while $\tau$ is the usual satellite twisting number.  If the local diagram is changed by a Reidemeister I move, then $\operatorname{wr}(D_K)$ changes by $\pm1$ and the box parameter must be changed by the compensating $\mp2$ half-twists in order to keep the same satellite twisting number.  In this paper the blackboard-framed local model $D_K$ is fixed once and for all; the notation $W(K,m)$ always refers to the direct limit built from that fixed model.

For $l\in \mathbb{Z}_{>0}$, let $T_l$ denote a solid torus standardly embedded (i.e. unknotted) in $S^3$. Given a solid torus $T'_l\subset T_l$, we aim to  construct an embedding $$h^{l+1}_{l}: T_{l}\rightarrow T_{l+1}.$$ To do this, we  denote  by $\alpha_l$ and $\beta_l$ the longitude and meridian  of $T_l$ and let $\delta_l$ and $\gamma_l$ be the longitude and meridian of $T'_l$ as illustrated in Figure \ref{3_1knot}. We then define the embedding $h^{l+1}_l: T_{l}\rightarrow T'_{l+1}\subset T_{l+1}$ so that $T_l$ is carried onto $T'_{l+1}$ with
\[h^{l+1}_l(\alpha_l)=\delta_{l+1}, \quad h^{l+1}_l(\beta_l)=\gamma_{l+1}.\]

 By the framing convention above, the embedding $h^{l+1}_l$ is chosen so that the box labelled $m$ is measured relative to the blackboard framing of the fixed diagram $D_K$ in the cube $C_l$.  Equivalently, the induced Whitehead double has satellite twisting number $\tau=m/2+\operatorname{wr}(D_K)$; see Proposition~\ref{infinity-surjection}.

 For any $l'>l$, we define the induced embedding $h^{l'}_{l}: T_l\rightarrow T_{l'}$
\begin{equation*}
  h^{l'}_{l}=h^{l'}_{l'-1}\circ h^{l'-1}_{l'-2}\circ \cdots\circ h^{l+2}_{l+1}\circ h^{l+1}_{l}.
\end{equation*}
Thus,  the collection  $(T_{l}, h^{l+1}_{l})$ forms a direct system, which admits  a direct limit.
The direct limit $W$ is defined as  the quotient space
\begin{equation*}
  q:  \sqcup_l T_l\rightarrow W,
\end{equation*}
where $q$ is the quotient map induced by the relation $\sim$ on $\sqcup_l T_l$. Explicitly, for $x\in T_s$ and $y\in T_{t}$, we set $x\sim y$ if and only if there exists an integer $k\geq \max\{s, t\}$ such that  $h^{k}_s(x)=h^{k}_{t}(y)$. Let $\iota_s: T_s\rightarrow \sqcup_l T_l$ denote the natural inclusion map. The composition $$q\circ \iota_s: T_s\rightarrow W$$
is an embedding. Consequently,  the images $q\circ\iota_l(T_l)$ form an exhaustion of $W$ by compact subsets, and therefore, $W$ is an open $3$-manifold. Let $T_l^*$ denote $q \circ \iota_l(T_l)$. $T_{l-1}^*$ is embedded in $T_{l}^*$ just as the way $h_{l-1}^{l}(T_{l-1})$ ($= T_{l}'$) is embedded in $T_{l}$.
 Hence, Figure \ref{3_1knot} can be viewed as a picture of the embedding of $T_{l-1}^*$ in $T_{l}^{*}$.

\subsection{Construction of contractible $3$-manifolds}
To obtain a contractible open 3-manifold, we will make use of specific knots inside the solid torus when constructing the embedding
$$h^{l+1}_{l}: T_{l}\rightarrow T_{l+1}.$$

Let $K$ be a nontrivial knot and  $T'_l$  be a tubular neighborhood of knot obtained by taking the connected sum of the Whitehead pattern with $m$ half-twists applied to the trivial knot with $K$,  as in Figure \ref{3_1knot}. As described in Section 3.1, this choice determines a family of embedding $h^{l+1}_{l}$. We then define
 $$W(K, m):=\lim_{l\to \infty}(T_l, h^{l+1}_{l}),$$ where
 $m$ is an even integer. By \cite[Proposition 2.1]{Gu21}, this construction yields a  contractible manifold.

\begin{proposition}
  For any $K$ and $m\in 2\mathbb{Z}$, the manifold $W(K, m
  )$ is contractible.
\end{proposition}

The topological structure of $W(K, m)$ is deeply intertwined with its fundamental group at infinity and the knot $K$. To elucidate their relationship, we consider  the image $T^*_l\subset W(K, m)$ induced by the composition $q\circ\iota_l : T_l \rightarrow W$. This construction yields a decomposition of $W(K, m)$
into amalgamation of $L_l = T^*_{l}\setminus T^*_{l-1}$'s. That is, for $l\geq 1$,
\begin{equation}\label{amalgamation of W}
W(K,m) = \lim_{l\to \infty} T_{0}^{\ast} \cup L_{1} \cup_{h_{1}^{2}} \cdots \cup_{h_{l-2}^{l-1}} L_{l-2} \cup_{h_{l-1}^{l}} L_l,\end{equation}
 where the sewing homeomorphism $h_{l}^{l+1}$ identifies the boundary component $\partial T_l$ of $L_l$ to the boundary component $\partial T'_{l+1}$of $L_{l+1}$.

This decomposition makes a bridge between topological information of the end of $W(K, m)$ and some knot groups.

\begin{proposition}\label{infinity-surjection} Let $\{K_l\}$ be a family of knots in $S^3$ with the following properties:
\begin{equation*}
  K_1=\mathcal{TW}_{\frac{m}{2}}\#K, \quad K_l=\operatorname{WD}_\tau(K_{l-1})\#K,
\end{equation*} where $\mathcal{TW}_{\frac{m}{2}}$
is an $\frac{m}{2}$-twist knot. Then, for any $l\geq 1$, there is a surjection
\begin{equation}\label{surjection}\pi_1(W(K, m)\setminus T^*_0)\rightarrow \pi_1(S^3\setminus K_l),\end{equation} where the satellite twisting number $\tau$ and the blackboard-framed half-twist parameter $m$ are related by
\begin{equation*}
  \tau=m/2+ \operatorname{wr}(D_K).
\end{equation*}
\end{proposition}

\begin{proof} Recall that for $0\leq s\leq l-1$, the embedding  $h_{l-s}^l: T_{l-s} \rightarrow T_l\subset S^3$ maps the solid torus $T_{l-s}$ into $T_l\subset S^3$.  The \emph{core} of a solid torus is defined as  a simple closed curve that serves as  a deformation retract of the solid torus.

We will proceed  by  induction on $s$ to show  that the core of the embedded solid torus  $h^{l}_{l-s}(T_{l-s})$ is isotopic to $K_s$ in $S^3$. This fact will be essential in constructing surjections between relevant fundamental groups.

For the base case $s=1$, the embedding $h^{l}_{l-1}:T_{l-1}\rightarrow T_{l}$ is illustrated in Figure \ref{3_1knot}. By construction, the core of $h^{l}_{l-1}(T_{l-1})$ is isotopic to the connected sum of the $\frac{m}{2}$-twist knot with $K$. Now assume that for some $s\geq 1$ the core of $h^{l}_{l-s}(T_{l-s})$ is isotopic to $K_s$. Consider the composition
\begin{equation*}
  h^{l}_{l-s-1}=h^l_{l-s}\circ h^{l-s}_{l-s-1}: T_{l-s-1}\rightarrow T_{l-s}\rightarrow T_l\subset S^3.
\end{equation*}
From Figure ~\ref{3_1knot},  the core of  $h_{l-s-1}^{l-s}(T_{l-s-1})$ is isotopic to $\mathcal{TW}_{\frac{m}{2}}\#K$.  Since $h^l_{l-s}(T_{l-s})$ is a tubular neighborhood of $K_s$, the image of $\mathcal{TW}_{\frac{m}{2}}$ under $h_{l-s}^l$ corresponds exactly to a $\tau$-twisted Whitehead double  of $K_s$, denoted by $\operatorname{WD}_\tau(K_s)$.  More precisely, the satellite twisting number $\tau$ is determined, with respect to the blackboard framing of $D_K$, by the Calugareanu--White--Fuller formula \cite{Whi69}
\begin{equation*}
  \tau = m/2 + \operatorname{wr}(D_K).
\end{equation*} Because the image of $K$ under $h^l_{l-s}$ remains isotopic to  $K$, we deduce that  the core of the image of  $h_{l-s-1}^l$ is isotopic to $K_{s+1}=\operatorname{WD}_\tau(K_s)\#K$, completing the induction.

We now use  the decomposition \eqref{amalgamation of W} to construct the required surjections.  From \eqref{amalgamation of W}, the fundamental group of the complement of $T^*_0$ admits a decomposition  as a free product with amalgamation: \begin{equation*}\pi_1(W(K, m)\setminus T^*_0)\cong\pi_1(T^*_1\setminus T^*_0)\ast_{\Lambda_1}\pi_1(T^*_2\setminus T^*_1)\ast_{\Lambda_2}\cdots\ast_{\Lambda_{l-1}}\pi_1(T^*_{l}\setminus T^*_{l-1})\ast_{\Lambda_l}\cdots,\end{equation*}
where $\Lambda_l=\pi_1(\partial T^*_l)$. For each $l\geq 1$, we consider the  natural quotient homomorphism:
\begin{equation*}
  \pi_1(W(K, m)\setminus T^*_0)\rightarrow \pi_1({W(K, m)}\setminus T^*_0)/ \textbf{N}_l,
\end{equation*} where $\textbf{N}_l$ is the normal subgroup generated by $\{\pi_1(T^*_{k}\setminus T^*_{k-1})\}_{k\geq l+1}$. By Claim 2 of \cite{Gu21}, there is a canonical isomorphism
\begin{equation*}
  \pi_1(W(K, m)\setminus T^*_0)/\textbf{N}_l\cong \pi_1(S^3\setminus K_l),
\end{equation*} and  hence we obtain a natural surjection $\pi_1(W(K,m)\setminus T^*_0)\rightarrow \pi_1(S^3\setminus K_l)$.
\end{proof}

\subsection{Proof of Theorem \ref{Thm: nonembeddability}}To apply Proposition \ref{infinity-surjection} effectively and to control the behavior of loops when passing to a large space, we require an additional local topological condition.

\begin{definition}
A topological space $X$ is \emph{locally \emph{1}-connected at the point} $x\in X$ if for each neighborhood $U$ of $x$,  there is a neighborhood $V\subset U$ of $x$ such that every loop in $V$ contracts in $U$. We say that $X$ is \emph{locally \emph{1}-connected} if $X$ is locally $1$-connected at each of its points.
\end{definition}
It was pointed out by  \cite[Lemma 1.1, P. 7]{Ste77} that such compact manifolds are characterized by their fundamental group.
\begin{lemma}\label{Lemma: Compact space f.g. pi_1}
  If $X$ is a compact, connected, locally connected, locally $1$-connected metric space, then $\pi_1(X)$ is finitely generated.
\end{lemma}

We combine this topological characterization with Proposition \ref{infinity-surjection} to complete the proof of Theorem \ref{Thm: nonembeddability}.

\begin{proof}[Proof of Theorem \ref{Thm: nonembeddability}] Let $K$ be  a nontrivial knot and $m$ be an even integer. Recall from Section~3.2 that  the contractible  $3$-manifold $W(K, m)$ is defined.  Suppose,  for contradiction that, $W(K, m)$ embedded into a compact, locally connected, locally $1$-connected metric space $X$.

Since the open solid torus $\text{Int}~T^*_0\subset W(K,m)$ is  pre-compact in $X$ where $T^*_0$ is constructed in Section~3.2, the complement $X\setminus \text{Int}~T^*_0$ is compact, locally connected, locally $1$-connected. By Lemma \ref{Lemma: Compact space f.g. pi_1}, $\pi_1(X\setminus \text{Int}~T^*_0)$ is finitely generated. Furthermore, by \cite[Claim 1]{Gu21}, there exist surjective homomorphisms: for any $l>0$
\[\pi_1(X\setminus \text{Int}~T^*_0)\twoheadrightarrow \pi_{1}(S^3\setminus K_l),\] where $K_l$ is defined in Proposition \ref{infinity-surjection}. In particular, the groups  $\pi_1(S^3\setminus K_l)$ have uniformly bounded rank.

On the other hand, we claim that for $l\geq 1$, the rank of $\pi_1(S^3 \setminus K_l)$
\begin{equation}\label{rank-bound}
  r(K_l)\geq l+1,
\end{equation} which contradicts the uniform bound above. Lemma \ref{Lemma: lower bound for the rank} reduces \eqref{rank-bound} to showing that the JSJ decomposition of  $S^3\setminus K_l$ has at least $l$ hyperbolic components.

We argue by induction on $l$. The case $l=0$ is trivial.  Assume that the statement holds for $l$. Recall  that the embedding $h^{l+1}_0$ factors as:
\begin{equation*}
h^{l+1}_0=h^{l+1}_{1}\circ h^1_0: T_0\rightarrow T_1\rightarrow T_{l+1},
\end{equation*} and hence
\begin{equation*}
  S^3\setminus h^{l+1}_0(T_0)=\left(S^3\setminus h^{l+1}_1(T_1)\right)\cup h^{l+1}_1\left(T_1\setminus h^1_0(T_0)\right).
\end{equation*}
From the proof of Proposition \ref{infinity-surjection}, the image of $h^{l+1}_s:T_s\rightarrow T_{l+1}$ is the tubular neighborhood of the knot $K_{l+s}$ for $s=0,1$ and thus
\begin{equation*}
  S^3\setminus h^{l+1}_0(T_0)\cong S^3\setminus K_{l+1}, \quad  S^3\setminus h^{l+1}_1(T_1)\cong S^3\setminus K_{l}.
\end{equation*}
It remains to analyze  $T_1\setminus h^{1}_0(T_0)$. From Figure~\ref{3_1knot}, let $\mathbf{B}$ denote the box containing the $K$-portion. Then
\begin{equation*}
  T_1\setminus h^{1}_0(T_0)=\left(T_1\setminus (h^{1}_0(T_0)\cup \mathbf{B})\right)\cup \left(\mathbf{B}\setminus h^{1}_0(T_0)\right).
\end{equation*}
The union $\mathbf{B}\cup h^{1}_0(T_0)$ is a tubular neighborhood of the Whitehead pattern with $m$ half-twists applied to the trivial knot, so $T_1\setminus (h^{1}_0(T_0)\cup \mathbf{B})$ is the complement of the Whitehead link with $m$ half-twists, which is  hyperbolic by Lemma \ref{Lemma: toral differences are hyperbolic}. Moreover, $\operatorname{Int} \mathbf{B}\setminus h^{1}_0(T_0)$ is homeomorphic to $S^3\setminus K$, and adding the incompressible tori in the JSJ decomposition of $S^3\setminus K$ yields a JSJ decomposition for $T_1\setminus h^{1}_0(T_0)$ with at least one hyperbolic part.

By inductive hypothesis,  $S^3\setminus h^{l+1}_1(T_1)$ contributes  at least $l$ hyperbolic parts and the above analysis yields one more. Thus, $S^3\setminus K_{l+1}$ has a JSJ decomposition with at least $l+1$ hyperbolic parts, completing the induction. Equation \eqref{rank-bound} follows, giving the desired contradiction.
\end{proof}

\begin{remark}
  Alternatively, we sketch an argument by combinatorial group theory. Since the rank of a group is at least as large as that of any homomorphic image, it suffices to show that the lower bound of the rank of $\pi_1(S^3 \setminus K_l)$ increases without bound as $l\to \infty$.
 By the Seifert--van Kampen theorem,
 \begin{equation*}\begin{aligned}
  G(\operatorname{WD}_{\tau}(K_2))
 & = G(\text{Whitehead link}) \ast_{\Lambda} G(\operatorname{WD}_{\tau}(K_1)\# K)\\&
  = G(\text{Whitehead link}) \ast_{\Lambda} G(\operatorname{WD}_{\tau}(\mathcal{TW}_{\frac{m}{2}} \# K)\#K),\end{aligned}
 \end{equation*}
 the amalgamated free product of the link group $G(\text{Whitehead link})$ and the knot group $G(\operatorname{WD}_{\tau}(K_1) \# K)$, with  the peripheral subgroup $\Lambda \cong \mathbb{Z}\oplus \mathbb{Z}$  corresponding to the torus along which spaces are glued via the homeomorphism $h$. For clarity, we may think of $G(\operatorname{WD}_{\tau}(K_2))$, $G(\text{Whitehead link})$, $G(\operatorname{WD}_{\tau}(K_1) \# K)$ and $\Lambda$ as $\pi_1(S^3 \setminus K_W)$, $\pi_1(V_P\setminus K_P)$ and $\pi_1(S^3 \setminus V_C)$ and $\pi_1(\partial V_P)$, respectively, all as described in Definition \ref{Def: Whitehead doubling}.

 Since the peripheral subgroup $\Lambda < G(\operatorname{WD}_{\tau}(K_1) \# K)$ is generated by the longitudinal and meridional generators, we can "abelianize" the knot group of $K$ to the meridional generator in $G(\operatorname{WD}_{\tau}(K_1 )\# K)$. By the universal property of amalgamated products, this induces a surjective homomorphism
\begin{equation*}\begin{aligned}
\phi: G(\operatorname{WD}_{\tau}(K_2)) =
G(\text{Whitehead link}) \ast_{\Lambda} &G(\operatorname{WD}_{\tau}(K_1) \# K) \\&\twoheadrightarrow G(\text{Whitehead link}) \ast_{\Lambda} G(\operatorname{WD}_{\tau}(K_1)).\end{aligned}\end{equation*}
 The image of $\phi$ is the knot group of a Whitehead double of $\operatorname{WD}_{\tau}(K_1)$,  denoted  by $G(\operatorname{WD}_{\tau}^{2}(K_1)):=G(\operatorname{WD}_{\tau}(\operatorname{WD}_{\tau}(K_1)))$. By iteration, $G(\operatorname{WD}_{\tau}^{l-1}(K_1))$ is a homomorphic image of $G(\operatorname{WD}_{\tau}(K_{l-1}))$, where $l\geq 3$. We observe that $G(\operatorname{WD}_{\tau}(K_{l-1})\#K)$ surjects onto  $G(\operatorname{WD}_{\tau}(K_{l-1}))$ via the "abelianization" of $G(K)$ onto the meridional generator. Combining this observation with  Proposition \ref{Thm: Lower bound for WD knot group} yields
 $$r(K_l) = r(\operatorname{WD}_{\tau}(K_{l-1})\# K) \geq r(\operatorname{WD}_{\tau}(K_{l-1}))\geq r(\operatorname{WD}_{\tau}^{l-1}(K_1))\geq l-1.$$
\end{remark}

\subsection{Non-homeomorphism}\label{Section: Non-homeomorphism}

The preceding argument shows that the topology of $W(K, m)$ encodes rich information about the knot and the blackboard-framed half-twist parameter $m$.  In particular, the marked JSJ structure at infinity forces strong rigidity properties: distinct choices of $(K,m)$ produce different marked layer structures.  This rigidity will be crucial in distinguishing contractible $3$-manifolds obtained from different pairs $(K,m)$. We establish a classification result asserting that the pair $(K,m)$ is a complete invariant for the homeomorphism type of $W(K, m)$.
\begin{theorem} \label{Thm: homeomorphism of manifolds}
Let $K$ and $K'$ be two nontrivial knots and let $m$ and $m'$ be even integers, with $W(K,m)$ and $W(K',m')$ constructed using the framing convention of Section~\ref{Section: nonembeddability}. Then $W(K, m)$ is homeomorphic to $W(K',m')$ if and only if
\begin{equation*}
  m=m' \quad \text{ and }\quad  K \text{ is isotopic to } K'.
\end{equation*}
\end{theorem}

To establish the theorem, we first observe that the choice of the knot $K$ determines the JSJ decomposition of $T^*_l\setminus T^*_0$, and this decomposition provides a way to distinguish different knots (see Theorem \ref{knot-choice-distinguish}). We then analyze the influence of the blackboard-framed half-twist parameter $m$ on the marked Whitehead-pattern component.  The relevant invariant is not just the hyperbolic component as an unmarked manifold, but the component together with the preferred longitude on its outer boundary and the regular-fiber slope inherited from the adjacent connected-sum Seifert component (see Theorem \ref{twisting-number}).

\begin{figure}[h!]
  \centering
  \includegraphics[ width=8cm, height=10cm]{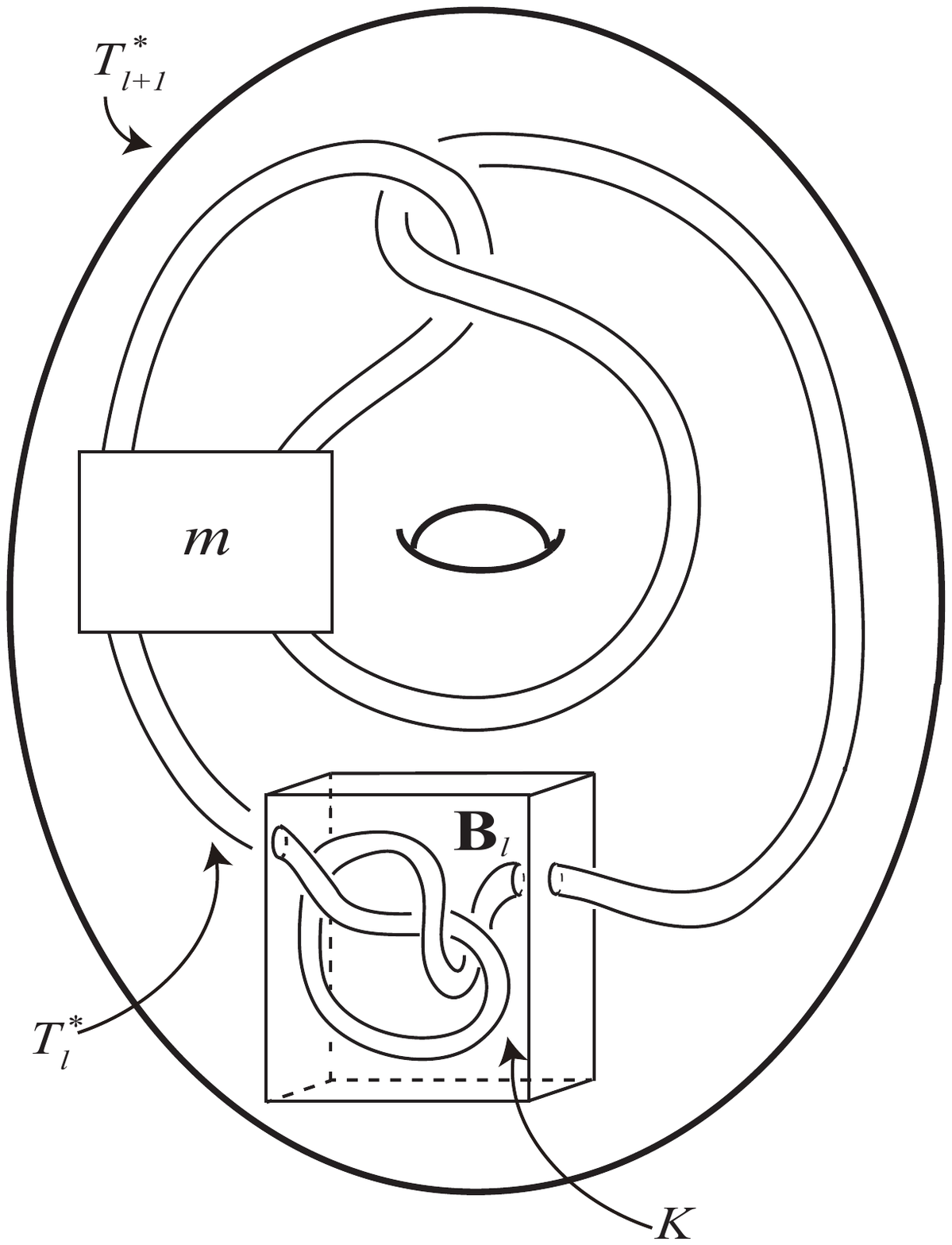}
  \caption{This picture illustrates the embedding of $T_{l}^{\ast}$ relative to $T_{l+1}^{\ast}$. The $\epsilon$-neighborhoods of $T_{l}^{\ast}$ and $\mathbf{B}_l$ are not shown.}
  \label{JSJ}
\end{figure}

\subsubsection{The knot choice }
Recall that for any knot $K$ and even integer $m$, the contractible manifold $W(K,m)$ admits a  decomposition
\begin{equation*}
  W(K,m)=T^*_0\cup (T^*_1\setminus T^*_0)\cup \cdots \cup (T^*_{l+1}\setminus T^*_l)\cup \cdots ,
\end{equation*} where $\{T^*_l\}_{l\geq 0}$ is an exhaustion of $W(K, m)$ by solid tori, as constructed in Section 3.2.  The proof of Proposition \ref{Thm: Lower bound for WD knot group} shows that  $T^*_{l+1}\setminus T^*_l$ contains a hyperbolic piece.

\begin{proposition}\label{existence-marked-hyperbolic} Let $K$ be a nontrivial knot. For any $l>0$, there exists a connected subset $\mathcal{H}^K_l\subset T^*_{l+1}\setminus T^*_l$ such that $\operatorname{Int}(\mathcal{H}^K_l)$ is hyperbolic, $\partial T^*_{l+1}\subset \mathcal{H}^K_l$ and
\begin{equation}\label{hyperbolic-isom}
  \text{the induced map }  H_2(\mathcal{H}_l^K)\rightarrow H_2(T^*_{l+1}\setminus T^*_l) \text{  is an isomorphism.}
\end{equation}
\end{proposition}
\begin{proof} Recall that for any $l\geq 0$, the solid torus $T^*_l$ is embedded into $T^*_{l+1}$ as illustrated in Figure \ref{JSJ}. In particular, the core of $T^{*}_l$ is the connected sum of the Whitehead-pattern piece with $m$ half-twists and $K$.

Let $\mathbf{B}_l$ denote the ball containing the part corresponding $K$. We may assume that $\partial \mathbf{B}_l$ intersects the core of $T^*_l$ transversely. Then, the intersection  $\partial \mathbf{B}_l\cap B(T^*_l, \epsilon) $ is a disjoint union of two discs, where $B(T^*_l, \epsilon)$ denotes the $\epsilon$-neighborhood of $T_{l}^{\ast}$. Observe that the union $\mathbf{B}_l\cup B(T^*_l,\epsilon)$ is a tubular neighborhood of the Whitehead-pattern piece with $m$ half-twists. We define that
\begin{equation*}
  \mathcal{H}^K_{l}:= T^*_{l+1}\setminus (B(\mathbf{B}_l, \epsilon)\cup B(T^*_l, \epsilon)).
\end{equation*}
It includes $\partial T^*_{l+1}$ and
its interior is the complement of the Whitehead link with $m$ half-twists, which is hyperbolic by Lemma \ref{Lemma: toral differences are hyperbolic}.

By the Mayer-Vietoris sequence, the map $H_2(\partial T^*_{l+1})\rightarrow H_2(T^*_{l+1}\setminus T^*_l)$ is an isomorphism. Since $\operatorname{Int} \mathcal{H}^K_l$ is the complement of the Whitehead link with $m$ half-twists, the analogous computation gives an isomorphism $H_2(\partial T^*_{l+1})\rightarrow H_2(\mathcal{H}^K_l)$. Therefore, we can conclude that the induced map $  H_2(\mathcal{H}_l^K)\rightarrow H_2(T^*_{l+1}\setminus T^*_l)$ is an isomorphism.
\end{proof}

By construction, the hyperbolic component of $T^*_{l+1}\setminus T^*_l$ comes from $\mathcal{H}^K_l$ and the JSJ decomposition of $S^3\setminus K$. However, $\mathcal{H}^K_l$ is the unique hyperbolic component with \eqref{hyperbolic-isom}.

\begin{corollary}\label{uniqueness-marked-hyperbolic} For any $l$, $\mathcal{H}^K_l$ is the unique hyperbolic component in $T^{*}_{l+1}\setminus T^*_l$ with the following property:
\begin{equation*}
\text{ the induced map }H_2(\mathcal{H}^K_l)\rightarrow H_2(T^*_{l+1}\setminus T^*_l) \text{ is an isomorphism. }
\end{equation*}
\end{corollary}
\begin{proof} Consider the decomposition of the space $T^*_{l+1}\setminus T^*_l$ as follows
\begin{equation}\label{JSJ-pre}
  T^*_{l+1}\setminus T^*_{l}= \mathcal{H}^K_l\cup \mathcal{S}^K_l\cup \mathbf{B}_l\setminus B(T^*_l, \epsilon),
 \end{equation} where  $\mathcal{S}^K_l$ denotes $B(T^*_l, \epsilon)\cup B(\mathbf{B}_l, \epsilon)\setminus (T^*_l\cup \mathbf{B}_l\setminus B(T^*_l, \epsilon))$. Observe that the interior $\mathcal{S}^K_l$ is homeomorphic to $S^1\times (S^2\setminus \{P_1, P_2, P_3\})$, which is Seifert fibered and the interior of $\mathbf{B}_l\setminus B(T^*_l, \epsilon)$ is homeomorphic to the complement of $K$ in $S^3$.

Consider a hyperbolic component $\mathcal{H}\subset T^*_{l+1}\setminus T^*_l$ which is not ambient isotopic to $\mathcal{H}^K_l$. From the above decomposition, we find that any hyperbolic component comes from $\mathcal{H}^K_l$ or the JSJ decomposition of $S^3\setminus K$. Thus, we have that $\mathcal{H}$ is ambient isotopic to some subset of $\mathbf{B}_l\setminus T^*_l$. Note that $H_2(\mathbf{B}_l\setminus B(T^*_l, \epsilon))\cong H_2(S^3\setminus K)\cong \{0\}$ and the induced map
\begin{equation*}
H_2(\mathcal{H})\rightarrow H_2(\mathbf{B}_l\setminus B(T^*_l, \epsilon))\rightarrow H_2(T^*_{l+1}\setminus T^*_l)
\end{equation*} is a trivial map, which completes the proof.
\end{proof}

\begin{definition}\label{marked-hyper} The subset $\mathcal{H}$ of $W(K, m)\setminus T^*_0$ is called the \emph{marked hyperbolic component} if it satisfies that
\begin{itemize}
  \item[(i)] $\partial \mathcal{H}$ has two components and its interior of $\mathcal{H}$ is hyperbolic; and
  \item[(ii)] the induced map $H_2(\mathcal{H})\rightarrow H_2(W(K, m)\setminus T^*_0)$ is an isomorphism.
\end{itemize}
\end{definition}

\begin{lemma}\label{marked-isotopic} Let $K$ be a nontrivial knot and let $\mathcal{H}$ be a marked hyperbolic component in $T^*_{l+1}\setminus T^*_0$. Then $\mathcal{H}$ is ambient isotopic to $\mathcal{H}^K_{l'}$ for some $l'\leq l$.
\end{lemma}
\begin{proof}
By the proof of Corollary \ref{uniqueness-marked-hyperbolic}, each hyperbolic component of $T^*_{k+1}\setminus T^*_{k}$ is ambient isotopic to $\mathcal{H}^K_{k}$ or a subset of $\mathbf{B}_k\setminus B(T^*_k, \epsilon)$. Combining this with the layer decomposition of  $T^*_{l+1}\setminus T^*_0$
\begin{equation*}
  T^*_{l+1}\setminus T^*_0=\bigcup^l_{k=0} (T^*_{k+1}\setminus T^*_{k})
\end{equation*} yields that
 any hyperbolic component of $T^*_{l+1}\setminus T^*_0$ is ambient isotopic to some $ \mathcal{H}^K_k$ or some subset of $\mathbf{B}_k\setminus B(T^*_k, \epsilon)$.

 Suppose that $\mathcal{H}$ is isotopic to  some subset of $\mathbf{B}_k\setminus B(T^*_k, \epsilon) $. Then the induced map
$$H_2(\mathcal{H})\rightarrow H_2(\mathbf{B}_k\setminus B(T^*_k,\epsilon))\rightarrow H_2(T^*_{k+1}\setminus T^*_k)$$ would be trivial, contradicting Definition \ref{marked-hyper}(ii), where $H_2(\mathbf{B}_k\setminus B(T^*_k,\epsilon))\cong H_2(S^3\setminus K)$ is trivial. Thus, $\mathcal{H}$ is isotopic to some $\mathcal{H}^K_k$.
\end{proof}

From Lemma \ref{marked-isotopic}, we know that one boundary component of any marked hyperbolic component is ambient isotopic to  some $\partial T^*_l$. We shall use this observation to distinguish among  knot choices.

\begin{theorem}\label{knot-choice-distinguish}
Let $K$ and $K'$ be nontrivial knots and let $m$ and $m'$ be even integers. If $W(K, m)$ is homeomorphic to $W(K', m')$, then one has that $K$ is ambient isotopic to $K'$.
\end{theorem}

\begin{proof} Since $W(K,m)\cong W(K', m')$, the space $W(K,m)$ admits another exhaustion $\hat{T}^*_l$ by solid tori and hence a second amalgamated decomposition
\begin{equation*}
  W(K, m)=\hat{T}^*_0\cup (\hat{T}^*_{1}\setminus \hat{T}^*_0)\cup \cdots \cup (\hat{T}^*_{l+1}\setminus \hat{T}^*_l)\cup\cdots,
\end{equation*} where the core of $\hat{T}^*_l$  is isotopic to the connected sum of the Whitehead-pattern piece with $m'$ half-twists and $K'$ in $\hat{T}^*_{l+1}$. Thus, there exist two integers $s,l>0$ such that
\begin{equation*}
  T^*_0\subset \hat{T}^*_{s}\subset \hat{T}^*_{s+1}\subset \hat{T}^*_{s+2}\subset {T}^*_{l+1}.
\end{equation*}
By Proposition \ref{existence-marked-hyperbolic}, there are two marked hyperbolic components, $\mathcal{H}^{K'}_s$ and $\mathcal{H}^{K'}_{s+1}$ of $T^*_{l+1}\setminus T^*_0$.

By Corollary \ref{uniqueness-marked-hyperbolic} and Lemma \ref{marked-isotopic}, $\mathcal{H}^{K'}_s$ is ambient isotopic to $\mathcal{H}^K_{l'}$ for some $l'\leq l$, implying that $\partial \hat{T}^*_{s+1}$ is ambient isotopic to $\partial {T}^*_{l'+1}$  in $T^*_{l+1}\setminus T^*_0$.\footnote{Orientation can be chosen to be preserved along these isotopies; in an orientable $3$-manifold, ambient isotopies of embedded tori can be taken to preserve the induced boundary orientation.} Similarly,
there exists an integer $l''>l'+1$ such that $\partial \hat{T}^*_{s+2}$ is ambient isotopic to $\partial T^*_{l''}$. Hence,
\begin{equation*}
  \hat{T}^*_{s+2}\setminus \hat{T}^*_{s+1}\cong T^*_{l''}\setminus T^*_{l'+1}.
\end{equation*}
 By Proposition \ref{existence-marked-hyperbolic}, the latter contains $l''-(l'+1)$ marked hyperbolic components in ${T}^*_{l''}\setminus {T}^*_{l'+1}$. The left-hand side contains exactly one marked hyperbolic component, namely $\mathcal{H}^{K'}_{s+1}$. Uniqueness (Corollary \ref{uniqueness-marked-hyperbolic}) forces that $l''-(l'+1)=1$, i.e. $l''=l'+2$. Thus,
\begin{equation*}
  \hat{T}^*_{s+2}\setminus \hat{T}^*_{s+1}\cong T^*_{l'+2}\setminus T^*_{l'+1}.
\end{equation*}
Combining  this correspondence with Corollary \ref{uniqueness-marked-hyperbolic}, we identify  $\mathcal{H}^K_{l'+1}$ with $\mathcal{H}^{K'}_{s+1}$. Thus,
\begin{equation}\label{key-homeomorphism}
  (\hat{T}^*_{s+2}\setminus \hat{T}^*_{s+1})\setminus \mathcal{H}^{K'}_{s+1}\cong (T^*_{l'+2}\setminus T^*_{l'+1})\setminus \mathcal{H}^K_{l'+1}.
\end{equation}
On the right-hand side, using \eqref{JSJ-pre} implies
\begin{equation}\label{JSJ-pre-again}
  (T^*_{l'+2}\setminus T^*_{l'+1})\setminus \mathcal{H}^K_{l'+1} \cong \mathcal{S}^K_{l'+1}\sqcup \mathbf{B}_{l'+1}\setminus B(T^*_{l'+1},\epsilon),
\end{equation}
where $\operatorname{Int}(\mathcal{S}^K_{l'+1})$ is $S^1\times (S^2\setminus \{P_1, P_2, P_3\})$.  Gluing a solid torus along $\partial T^*_{l'+1}\subset \partial \mathcal{S}^K_{l'+1} $ yields
\begin{equation}\label{product-topology}
  \operatorname{Int}~\mathcal{S}^K_{l'+1}\cup_{\partial T^*_{l'+1}} S^1\times D^2 \cong T^2\times (0, 1),
\end{equation} while the interior of $\mathbf{B}_{l'+1}\setminus B(T^*_{l'+1},\epsilon)$ is homeomorphic to $S^3\setminus K$. Consequently,
\begin{equation}\label{right-hand}
  (T^*_{l'+2}\setminus T^*_{l'+1})\setminus \mathcal{H}^K_{l'+1}\bigcup_{\partial T^*_{l'+1}} S^1\times D^2\cong S^3\setminus K.
\end{equation}
An analogous argument shows that
\begin{equation}\label{left-hand}
  (\hat{T}^*_{s+2}\setminus \hat{T}^*_{s+1})\setminus \mathcal{H}^{K'}_{s+1}\bigcup_{\partial \hat{T}^*_{s+1}} S^1\times D^2\cong S^3\setminus K'.
\end{equation}
Combining \eqref{right-hand} \eqref{left-hand} with \eqref{key-homeomorphism} yields
\begin{equation*}
  S^3\setminus K\cong S^3\setminus K'
\end{equation*} and  therefore  $K$ is isotopic to $K'$.
\end{proof}

\subsubsection{The twisting number}
We next recover the even integer $m$ from the marked Whitehead-pattern component together with its natural boundary slopes.  The point is that the unmarked hyperbolic piece is not the whole invariant: one must also remember how this piece is attached to the surrounding solid-torus and connected-sum data.

For a slope $\rho$ on a torus boundary component of a compact $3$-manifold $M$, we write $M(\rho)$ for the Dehn filling of that boundary component along $\rho$.  For the layer
\begin{equation*}
  L_l(K,m):=T^*_{l+1}\setminus \operatorname{Int}T^*_l,
\end{equation*}
we write $\mathcal{H}^{K,m}_l$ for the marked hyperbolic component previously denoted by $\mathcal{H}^K_l$, in order to display the dependence on $m$.  Put
\begin{equation*}
  O_l:=\partial T^*_{l+1}\subset \partial\mathcal{H}^{K,m}_l,\qquad
  Q_l:=\mathcal{H}^{K,m}_l\cap \mathcal{S}^K_l .
\end{equation*}
Thus $O_l$ is the outer boundary component of the marked Whitehead-pattern component, while $Q_l$ is the JSJ torus along which $\mathcal{H}^{K,m}_l$ meets the connected-sum Seifert component $\mathcal{S}^K_l$.

The slope on $O_l$ used below is not the meridian of the solid torus $T^*_{l+1}$, but its preferred longitude.  Equivalently, it is the slope along which one attaches the complementary solid torus of $T^*_{l+1}$ when $T^*_{l+1}$ is viewed as a standard unknotted solid torus in $S^3$.  The following lemma records the intrinsic form of this slope inside the open manifold.

\begin{lemma}\label{end-longitude-slope}
Let
\begin{equation*}
  E_l^+=W(K,m)\setminus \operatorname{Int}T^*_{l+1}.
\end{equation*}
Then the kernel of the map
\begin{equation*}
  H_1(O_l;\mathbb Z)\longrightarrow H_1(E_l^+;\mathbb Z)
\end{equation*}
induced by inclusion is generated by the preferred longitude of the solid torus $T^*_{l+1}$.  We denote this slope by $\lambda_l^+$ and call it the outer end-longitude of the marked component $\mathcal{H}^{K,m}_l$.
\end{lemma}

\begin{proof}
Let $\alpha_{l+1}$ and $\beta_{l+1}$ denote, respectively, the preferred longitude and meridian of $T^*_{l+1}$, using the notation of the construction.  Since the Whitehead pattern has algebraic winding number zero, the core of $T^*_{l+1}$ has algebraic winding number zero in $T^*_N$ for every $N>l+1$.  Hence $\alpha_{l+1}$ is null-homologous in $T^*_N\setminus \operatorname{Int}T^*_{l+1}$ for every sufficiently large $N$, and therefore maps trivially in $H_1(E_l^+;\mathbb Z)$.

Conversely, the Mayer--Vietoris sequence for
\begin{equation*}
  W(K,m)=T^*_{l+1}\cup E_l^+,\qquad T^*_{l+1}\cap E_l^+=O_l,
\end{equation*}
together with the contractibility of $W(K,m)$, shows that the kernel of $H_1(O_l;\mathbb Z)\to H_1(E_l^+;\mathbb Z)$ is a primitive rank-one subgroup whose image in $H_1(T^*_{l+1};\mathbb Z)$ is an isomorphism.  Thus this kernel is generated by a longitude of $T^*_{l+1}$.  Since $\alpha_{l+1}$ lies in the kernel, it generates the kernel.
\end{proof}

Let $\varphi_l$ denote the slope on $Q_l$ represented by the regular fiber of
\begin{equation*}
  \operatorname{Int}\mathcal{S}^K_l\cong S^1\times (S^2\setminus\{P_1,P_2,P_3\}).
\end{equation*}
The slope $\varphi_l$ is well-defined up to orientation, since the Seifert fibration of $S^1\times (S^2\setminus\{P_1,P_2,P_3\})$ is unique up to isotopy.

\begin{lemma}\label{slope-filling-twist}
With the notation above, there is a homeomorphism of marked exteriors
\begin{equation*}
  (\mathcal{H}^{K,m}_l(\lambda^+_l), Q_l, \varphi_l)
  \cong
  (S^3\setminus \operatorname{Int}N(\mathcal{TW}_{m/2}),
  \partial N(\mathcal{TW}_{m/2}), \mu_{\mathcal{TW}_{m/2}}),
\end{equation*}
where $\mathcal{TW}_{m/2}$ is the $m/2$-twist knot and $\mu_{\mathcal{TW}_{m/2}}$ is its meridian slope.
\end{lemma}

\begin{proof}
By construction, before the $K$-summand is attached, $\mathcal{H}^{K,m}_l$ is the exterior in the unknotted solid torus $T^*_{l+1}$ of the $m/2$-twist-knot summand.  Filling the boundary component $O_l=\partial T^*_{l+1}$ along the preferred longitude $\lambda^+_l$ attaches the complementary solid torus of $T^*_{l+1}$ in $S^3$.  The resulting one-cusped manifold is therefore the exterior of the $m/2$-twist knot.

It remains to identify the marked slope on the remaining boundary.  The component $\mathcal{S}^K_l$ is the standard $P\times S^1$ piece appearing in the connected-sum decomposition of the exterior of $\mathcal{TW}_{m/2}\# K$, where $P$ is a pair of pants.  In this decomposition the regular fiber is the meridian slope on the boundary of each summand exterior.  Hence the slope $\varphi_l$ on $Q_l$ is identified, after the filling along $\lambda^+_l$, with the meridian $\mu_{\mathcal{TW}_{m/2}}$ of the twist knot.
\end{proof}

\begin{lemma}\label{slope-preservation}
Let a homeomorphism between two such manifolds be isotoped so that it sends a layer $L_l(K,m)$ to a layer $L_s(K',m')$ and sends $\mathcal{H}^{K,m}_l$ to $\mathcal{H}^{K',m'}_s$.  Then it sends the outer end-longitude $\lambda^+_l$ to $\pm\lambda^{+\prime}_s$ and sends the fiber slope $\varphi_l$ to $\pm\varphi'_s$.
\end{lemma}

\begin{proof}
The slope $\lambda^+_l$ is characterized by Lemma~\ref{end-longitude-slope} as the primitive kernel of
\begin{equation*}
  H_1(O_l;\mathbb Z)\longrightarrow H_1(W(K,m)\setminus \operatorname{Int}T^*_{l+1};\mathbb Z).
\end{equation*}
A homeomorphism carrying the matched layer $L_l(K,m)$ to $L_s(K',m')$ carries the outside of $T^*_{l+1}$ to the outside of the corresponding torus in the target.  It therefore carries this kernel to the corresponding kernel, and hence sends $\lambda^+_l$ to $\pm\lambda^{+\prime}_s$.

The torus $Q_l$ is characterized as the JSJ torus through which the marked hyperbolic component is attached to the connected-sum Seifert component $\mathcal{S}^K_l$.  By uniqueness of the JSJ decomposition, the homeomorphism sends $\mathcal{S}^K_l$ to the corresponding component $\mathcal{S}^{K'}_s$.  Since this component is $S^1$ times a pair of pants, its Seifert fibers are intrinsic up to orientation.  Thus the regular-fiber slope on $Q_l$ is carried to the regular-fiber slope on $Q'_s$.
\end{proof}

\begin{theorem}\label{twisting-number}
  Let $K$ be a nontrivial knot and let $m$ and $m'$ be even integers. If $W(K, m)\cong W(K, m')$, then one has that $m=m'$.
\end{theorem}

\begin{proof}
Let $F:W(K,m')\rightarrow W(K,m)$ be a homeomorphism.  Applying the layer-matching argument in the proof of Theorem~\ref{knot-choice-distinguish} with $K'=K$, and then composing $F$ with an ambient isotopy of the target if necessary, we may assume that for some $s,l$ the image of a layer of the $m'$-exhaustion is a layer of the $m$-exhaustion and that
\begin{equation*}
  F(\widehat{\mathcal{H}}^{K,m'}_{s})=\mathcal{H}^{K,m}_{l},
\end{equation*}
where $\widehat{\mathcal{H}}^{K,m'}_{s}$ denotes the source marked component formed using the blackboard-framed half-twist parameter $m'$ and $\mathcal{H}^{K,m}_{l}$ denotes the target marked component formed using the blackboard-framed half-twist parameter $m$.  Moreover, the layer matching sends the outer boundary component to the outer boundary component and the connected-sum JSJ torus to the connected-sum JSJ torus.

By Lemma~\ref{slope-preservation}, $F$ sends the outer end-longitude of the source marked component to the outer end-longitude of the target marked component.  Hence $F$ extends over the corresponding Dehn fillings of the outer boundary components.  Lemma~\ref{slope-filling-twist} therefore gives a meridian-preserving homeomorphism
\begin{equation*}
  S^3\setminus \operatorname{Int}N(\mathcal{TW}_{m'/2})
  \cong
  S^3\setminus \operatorname{Int}N(\mathcal{TW}_{m/2}).
\end{equation*}
Consequently the Alexander polynomials of these two twist knots are equal, up to multiplication by a unit in $\mathbb{Z}[t,t^{-1}]$.

Let $n=m/2$ and $n'=m'/2$.  For the $n$-twist knot, Rolfsen's computation \cite[Ex. 7, P. 166]{Rol76} gives, in the symmetric normalization,
\begin{equation*}
  \Delta_{\mathcal{TW}_{n}}(t)=n(t+t^{-1})+1-2n
\end{equation*}
(equivalently, $\Delta_{\mathcal{TW}_{n}}(t)= n t^2+(1-2n)t+n$).  Equality of the normalized Alexander polynomials therefore implies $n=n'$, and hence $m=m'$.
\end{proof}

Theorem \ref{Thm: homeomorphism of manifolds} follows readily from Theorems \ref{knot-choice-distinguish} and \ref{twisting-number}.  Combining it with Theorem \ref{Thm: nonembeddability}, we obtain the following result.
\begin{corollary}\label{Cor: infinitely many W}
  There exist infinitely many non-homeomorphic contractible open $3$-manifolds that do not embed as an open subset in any compact, locally connected, locally $1$-connected metric $3$-space.
\end{corollary}

\begin{proof}[Proof of Theorem \ref{Thm: high dimensional nonembeddability}]
Utilizing the proof of \cite[Thm. 1.2]{Gu21}, we may extend Theorem \ref{Thm: nonembeddability} and Theorem \ref{twisting-number} in the same way that Theorem 1.2 extends Theorem 1.1 in \cite{Gu21}. Therefore, we obtain a generalization of Corollary \ref{Cor: infinitely many W}.
\end{proof}

\begin{figure}[h!]
  \centering
  \includegraphics[ width=7cm, height=9cm]{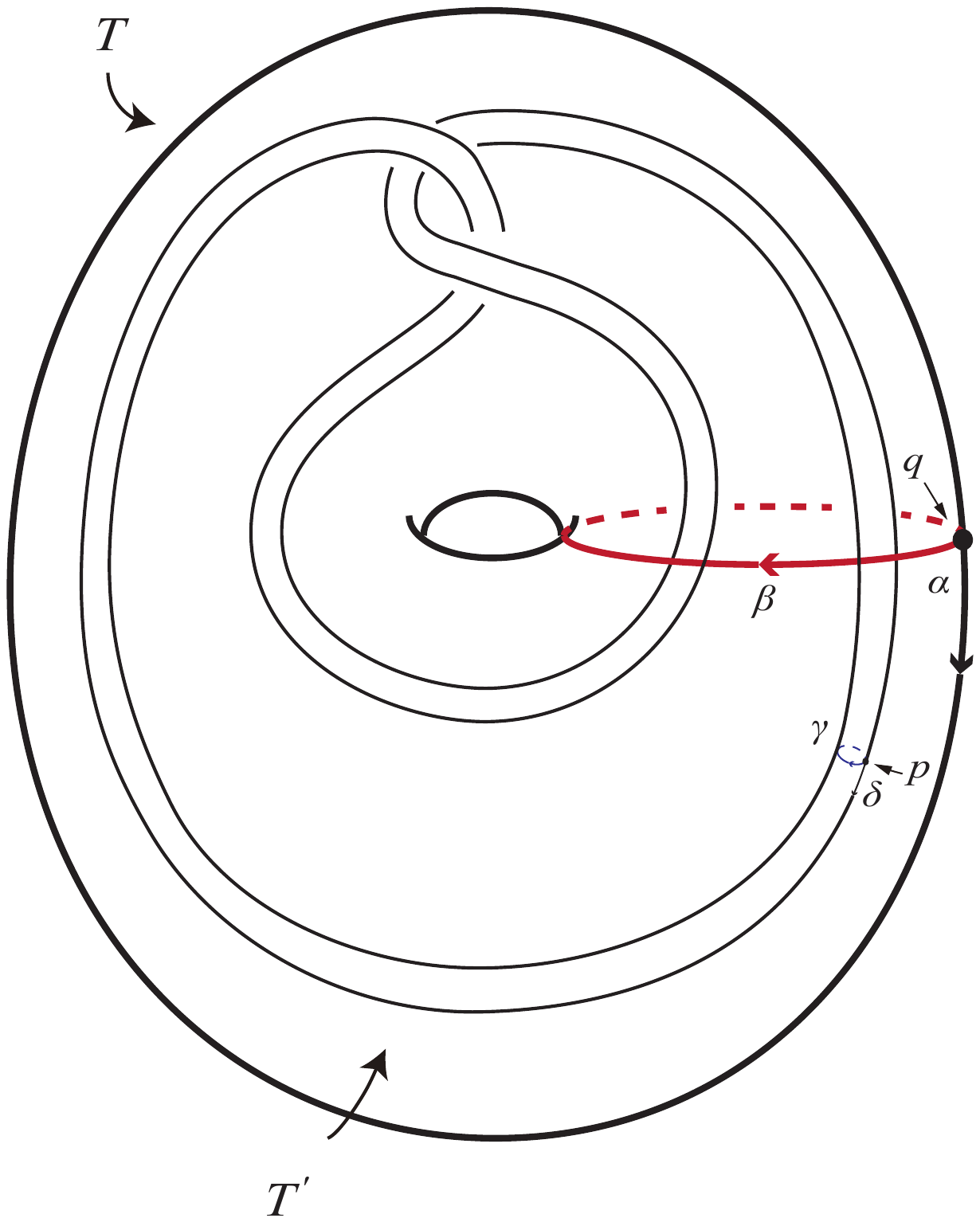}
  \caption{$L^\ast = T \setminus T'$. The ``inner'' boundary component of $L^\ast$ is $\partial T'$. The ``outer'' boundary component
  of $L^\ast$ is $\partial T$.}
  \label{whitehead}
\end{figure}

Recall the amalgamation of $W^3$ in (\ref{amalgamation of W}).
Setting the number of half-twists $m =0$ and unknotting the cube with a trefoil-knotted hole, as shown in Figure \ref{3_1knot} results in a cobordism $L^\ast$, which is widely known as the first stage of constructing a Whitehead manifold. See Figure \ref{whitehead}.

Consider a variation of $W^3$ by inserting infinitely many copies of $L^*$  in (\ref{amalgamation of W})
\begin{equation}\label{modified amalgamation of W}
W' = \lim_{j\to \infty} T_{0}^{\ast} \cup L_{1} \cup  \cdots \cup L_{l} \cup_{H_{\ast}^{l}}  L^\ast
\cup_{H_{l+1}^{\ast}} L_{l+1}\cup \cdots \cup L_{j}
\end{equation}
where the sewing homeomorphism $H_{*}^{l}$ identifies the boundary component $\partial T_l$ of $L_l$ to the boundary component $\partial T'$
of $L^\ast$ and the sewing homeomorphism  $H_{l+1}^{*}$ identifies the boundary component $\partial T$ of $L^\ast$ to the boundary component $\partial T'_{l+1}$ of $L_{l+1}$. Symbolically, let $L$ denote a copy of $L_l$, and let $+$ denote the union operator in (\ref{modified amalgamation of W}). When $L_l$ contains a trefoil-knotted hole (with $m=0$), a contractible open 3-manifold constructed by R. H. Bing can be written as
\begin{equation}\label{Bing's manifold}
\text{Bing's manifold}= T_{0}^{\ast} + L + L + L+\cdots
\end{equation}
A variation, which is a contractible open 3-manifold constructed by Sternfeld in \cite{Ste77}, can be written as
\begin{equation}\label{Sternfeld's manifold}
\text{Sternfeld's manifold} = T_{0}^{\ast} + L + L^\ast + L + L^\ast + \cdots
\end{equation}
 By inserting infinitely many copies of $L^*$ in (\ref{amalgamation of W}), we may obtain an infinite collection $\mathcal{C}$. For instance,
 $$W'' = T_{0}^{\ast} + L + L^\ast + L + L^\ast + L^\ast + L +\cdots$$
 It follows from \cite[Prop. 7]{Gu21} that
Sternfeld's manifold belongs to this collection $\mathcal{C}$. Because the nonembeddability of all manifolds in $\mathcal{C}$, as in the proof of Theorem \ref{Thm: nonembeddability}, can be reduced to estimating the lower bound of iterated Whitehead doubles of a nontrivial knot, we affirmatively answer \cite[Question 1]{Gu21}.
\begin{corollary}
  Each manifold in $\mathcal{C}$ embeds in no compact, locally connected and locally $1$-connected metric $3$-space.
\end{corollary}

\end{document}